\documentclass[a4paper]{article}
\usepackage[style=authoryear]{biblatex}
\usepackage{graphicx}%
\usepackage{multirow}%
\usepackage{amsmath,amssymb,amsfonts}%
\usepackage{amsthm}%
\usepackage{authblk}%
\usepackage{mathrsfs}%
\usepackage[title]{appendix}%
\usepackage{xcolor}%
\usepackage{textcomp}%
\usepackage{manyfoot}%
\usepackage{booktabs}%
\usepackage{algorithm}%
\usepackage{algorithmicx}%
\usepackage{algpseudocode}%
\usepackage{listings}%
\usepackage{hyperref}

\usepackage{subcaption}%

\theoremstyle{thmstyleone}%
\newtheorem{theorem}{Theorem}%

\begin{document}

\title{A hybrid optimization framework for the General Continuous Energy-Constrained Scheduling Problem}

\author[1]{Roel Brouwer}

\author[1]{Marjan van den Akker}

\author[1]{Han Hoogeveen}

\affil[1]{Department of Information and Computing Sciences, Utrecht University, Utrecht, The Netherlands}

\begin{abstract}
	We present a hybrid optimization framework for a class of problems, formalized as a generalization of the Continuous Energy-Con\-strained Scheduling Problem (CECSP), introduced by Nattaf et al. (2014). This class is obtained from challenges concerning demand response in energy networks. Our framework extends a previously developed approach. A set of jobs has to be processed on a continuous, shared resource. Consequently, a schedule for a job does not only contain a start and completion time, but also a resource consumption profile, where we have to respect lower and upper bounds on resource consumption during processing. In this work, we develop a hybrid approach for the case where the objective is a step-wise increasing function of completion time, using local search, linear programming and $O(n)$ lower bounds. We exploit that the costs are known in the local search and use bounds to assess feasibility more efficiently than by LP. We compare its performance to a mixed-integer linear program. After that, we extend this to a hybrid optimization framework for the General CECSP. This uses an event-based model, and applies a decomposition in two parts: 1) determining the order of events and 2) finding the event times, and hence the start and completion times of jobs, together with the resource consumption profiles. We argue the broad applicability of this framework.
\end{abstract}

\maketitle

\section{Introduction} \label{sec:gcecsp-intro}
In this paper we study scheduling problems related to demand response in energy networks. For example, when charging the battery of an e-vehicle, it must be provided with enough energy to fill it up to capacity, but the exact charging time can be freely determined as long as charging takes place between the arrival and departure time of the vehicle. Moreover, we are not necessarily limited to a fixed rate for charging the battery, although we want to take into account that preemption of charging may be undesirable, due to the response time of batteries.

This leads us to consider scheduling problems where we are primarily interested in the (cumulative) amount of work (or resource) required to complete a job, while the rate of consumption may vary within a given range. Such problems require a different approach than the more traditional ones, as both resources and time are continuously-divisible, and jobs are flexible (to a certain extent) along both axes. This type of problem, that we will refer to as the General Continuous Energy-Constrained Scheduling Problem (GCECSP), has not yet been studied extensively.

The GCECSP is described as follows. A set $\{J_1,\ldots,J_n\}$ of jobs has to be processed on a continuous, shared resource $R$. This means that, at any time, multiple jobs can be processed simultaneously and at different rates, as long as their total consumption does not exceed the available resource capacity $P$. A schedule for a job $J_j$ does not only contain a start and completion time (respecting the release time $r_j$ and deadline $\bar{d}_j$), but also a resource consumption profile $p_j(t)$, where we have to respect lower and upper bounds ($P^-_j,P^+_j$) on resource consumption during processing. The total consumption of job $J_j$ should equal its cumulative requirement $E_j$. Preemption is not allowed: from its start until its completion, each job must consume resources at a rate of at least $P^-_j > 0$ units.

The GCECSP is a generalization of the Continuous Energy-Constrained Scheduling Problem (CECSP) that we studied in our previous work \parencite{Brouwer2023} and was originally introduced by \cite{Nattaf2014}. The CECSP is a generalization of the Cumulative Scheduling Problem (CuSP) presented by \cite{Baptiste1999}. In the case of the CuSP, a number of activities have to be scheduled on a single shared resource with a given capacity. Each activity has a release time, deadline, (fixed) processing time and (constant) resource capacity requirement. The CuSP was formulated as a subproblem of the Resource Constrained Project Scheduling Problem (RCPSP), relaxing precedence constraints and considering only a single resource. The RCPSP is a very general problem that concerns the scheduling of activities subject to precedence, time and resource constraints. The surveys by \cite{Hartmann2010,Hartmann2022} provide a good overview of the RCPSP and its extensions. Most closely related to the present work, continuous and event-based formulations of the RCPSP have been studied and evaluated by \cite{Kone2011} and \cite{Kopanos2014} as well as more general models with flexible resource profiles in continuous time (FRCPSP), for example by \cite{Naber2017}. 

For the CECSP including efficiency functions that apply to the conversion of the amount of consumed resource to the contribution towards the resource requirement of a job, with minimizing resource consumption as the objective, \cite{Nattaf2017} and \cite{Nattaf2019} provided several algorithms aiming to find exact solutions for small instances. The CECSP generalizes the CuSP as the resource capacity requirement is considered to be a range with a lower and an upper bound, rather than a fixed value, and the consumption rate can vary during the execution of the job. As a result, the processing time can vary, depending on the consumption rate during execution. Through its relation to CuSP, \cite{Nattaf2016} proved that the problem of finding a feasible solution for CECSP is NP-complete.

The GCECSP is also closely related to the scheduling of malleable jobs (as introduced by \cite{Turek1992}) on parallel machines, which involves the scheduling of jobs on $P$ machines, while the number of machines assigned to a job can change during its execution. The machines can be viewed as a discretized resource.

A wide range of objectives can be considered in variants of the GCECSP. In this work, we explore the GCECSP with step-wise cost function in detail. Problems with a generalized step-wise cost function have not been widely studied. In scheduling literature, \cite{Detienne2011,Detienne2012} study single ($1||\bar{f}(C_j)$) and parallel ($R|r_j|\bar{f_j}(C_j)$) machine scheduling problems with a regular step-wise cost function. The single machine problem is proven to be NP-hard in the strong sense. For both the single machine and the parallel unrelated machine problem, \cite{Detienne2011} introduce a dominance relation, proving that if there exists at least one feasible schedule in which each job completes before its deadline, then the schedules where jobs are sequenced in non-decreasing order of their deadlines are feasible. \cite{Tseng2010} study a similar problem on a single machine, but refer to it as \textit{stepwise tardiness}.

\textbf{Our contribution} is twofold. (1) We present an optimization framework for the GCECSP and related problems, based on our previously developed hybrid local search approach \cite{Brouwer2023}, using simulated annealing and linear programming. (2) As an example, we detail how our approach can be applied to a variant with a step-wise constant objective function, where we exploit the properties of the objective to improve the efficiency of the approach.

The rest of this work is structured as follows. First, we will give a detailed problem description in Section \ref{sec:gcecsp-prob-desc}. Then, we will introduce our framework in Section \ref{sec:gcecsp-hybrid-alg}, using an interesting variant of the problem with step-wise cost functions as an example throughout. In this section, we will discuss the event-based model, a proposed decomposition of the problem, and the implementation of our solution approach. Our approach uses local search, linear programming and a number of algorithms for approximating the penalty terms, to avoid solving the LP in every iteration. We continue with an overview of variants of the GCECSP, and a discussion on the applicability of our framework to these variants in Section \ref{sec:gcecsp-variants}. In Section \ref{sec:gcecsp-results}, we will describe the generation of test instances for the problem introduced in Section \ref{sec:gcecsp-hybrid-alg} and present the computational results of applying our approach to them. Finally, we will draw conclusions in Section \ref{sec:gcecsp-conclusion}. An overview of the notation used in this work can be found in Table \ref{tab:gcecsp-notation}.
\begin{table}
\caption{Overview of notation} \label{tab:gcecsp-notation}
\centering
	\begin{tabular}{@{}ll@{}}
		\toprule
		\multicolumn{2}{l}{Indices} \\
		\midrule
		$j$ & job \\
		$i$ & event, interval \\
		$e$ & position in event order \\
		$k$ & jump point index \\
		$I(e)$ & converts position to index \\
		$E(i)$ & converts index to position \\
		\midrule
		\multicolumn{2}{l}{Sets} \\
		\midrule
		$J$ & Set of jobs \\
		$\mathcal{E}$ & List of events \\
		$\mathcal{F}$ & Fixed-time events \\
		$\mathcal{P}$ & Plannable events \\
		\midrule
		\multicolumn{2}{l}{Parameters} \\
		\midrule
		$n$ & number of jobs \\
		$k$ & number of cost intervals \\
		$P$ & available amount of resource $R$ \\
		$E_j$ & resource requirement of job $J_j$ \\
		$r_j, \bar{d}_j$ & release time and deadline of job $J_j$ \\
		$P^-_j, P^+_j$ & lower and upper bound on resource consumption for job $J_j$ \\
		$S_j, C_j$ & start and completion time of job $J_j$ \\
		$K^k_j$ & jump point $k$ of job $J_j$ \\
		$\bar{f}_j$ & step-wise increasing cost function of job $J_j$ \\
		$w_{j,k}$ & additional cost for completing job $J_j$ after jump point $K^{k-1}_j$ \\
		$p_j(t)$ & resource consumption profile of job $J_j$ \\
		$M$ & large constant term (`big M')\\
		$L^B, L^R$ & penalty for violating bounds (B) or resource availability (R) \\
		\midrule
		\multicolumn{2}{l}{Variables} \\
		\midrule
		$t_i$ & time of occurrence of the event with index $i$ \\
		$p_{j,i}$ & amount of resource consumed by job $J_j$ in the interval that starts when the \\
		& event with index $i$ occurs \\
		$a_{i,i'}$ & $\begin{cases} 1 & \text{if } E(i) < E(i') \\ 0 & \text{if } E(i) \geq E(i') \end{cases}$ relative order of events $i$ and $i'$\\
		$b_{i,i'}$ & $\begin{cases} 1 & \text{if } i' = \arg\min_{i' \in \mathcal{P}}(E(i') - E(i)) \\ 0 & \text{otherwise} \end{cases}$ $\begin{aligned}&\text{successor relation of plannable}\\ &\text{events }i\text{ and }i'\end{aligned}$\\
		$s^-_{j,i}, s^+_{j,i}$ & amount of resource job $J_j$ consumes less (more) than its lower (upper) \\
		& bound in interval $i$\\
		$s^t_i$ & amount of resource consumed in interval $i$ above the available amount $P$\\
		\bottomrule
	\end{tabular}
\end{table}

\section{Problem description} \label{sec:gcecsp-prob-desc}
We consider a wide range of scheduling problems. The common basis of these problems, which we will call the General Continuous Energy-Constrained Scheduling Problem (GCECSP) from here on, is defined by the following properties:
\begin{itemize}
	\item A resource $R$ with a constant availability of $P$;
	\item For each job $J_j,$ $j \in \{1, ..., n \}$:
	\begin{itemize}
		\item Resource requirement $E_j$;
		\item Release time $r_j$;
		\item Deadline $\bar{d}_j$;
		\item Lower bound $P^-_j$;
		\item Upper bound $P^+_j$.
	\end{itemize}
\end{itemize}

For each job, we need to determine a start time $S_j$, a completion time $C_j$ and a resource consumption profile $p_j(t)$ (a function that describes the amount of resource that job $J_j$ consumes over time) such that \textit{some objective} (most commonly a function of $C_j$) is minimized, while the following constraints are respected:
\begin{itemize}
	\item[C1] The total amount of resource job $J_j$ consumes is exactly $E_j$, i.e. $\int_t p_j(t)\ dt = E_j$;
	\item[C2] Job $J_j$ does not start before its release time $r_j$, i.e. $S_j \geq r_j$;
	\item[C3] Job $J_j$ completes no later than its deadline $\bar{d}_j$, i.e. $C_j \leq \bar{d}_j$;
	\item[C4] Job $J_j$ only consumes resources between its start and completion time, i.e. $p_j(t) = 0$ for $t < S_j$ or $t \geq C_j$;
	\item[C5] While active, the resource consumption of  job $J_j$ never drops below $P_j^-$ or rises above $P_j^+$, i.e. $P_j^- \leq p_j(t) \leq P_j^+$ for $S_j \leq t \leq C_j$;
	\item[C6] The total amount of resource consumed by all jobs together at any given time can never exceed $P$, i.e. $\sum_j p_j(t) \leq P$ for all possible values of $t$.
\end{itemize}

\begin{figure}
	\centering
	\tabskip=0pt
	\valign{#\cr
		\hbox{%
			\begin{subfigure}{.55\textwidth}
				\centering
				\begin{tabular}[b]{@{}l|lllll@{}}
					\toprule
					$j$ & $E_j$ & $r_j$ & $\bar{d}_j$ & $P^-_j$ & $P^+_j$ \\
					\midrule
					1 & 70.0 & 0.0 & 3.0 & 10.0 & 30.0 \\
					2 & 20.0 & 1.0 & 3.0 & 10.0 & 40.0 \\
					3 & 45.0 & 2.5 & 4.0 & 10.0 & 50.0 \\
					\bottomrule
				\end{tabular}
				\caption{Job properties}
				\label{fig:gcecsp-example-sub-properties}
			\end{subfigure}%
		}\vfill
		\hbox{%
			\begin{subfigure}{.55\textwidth}
				\centering
				\includegraphics[width=0.9\textwidth]{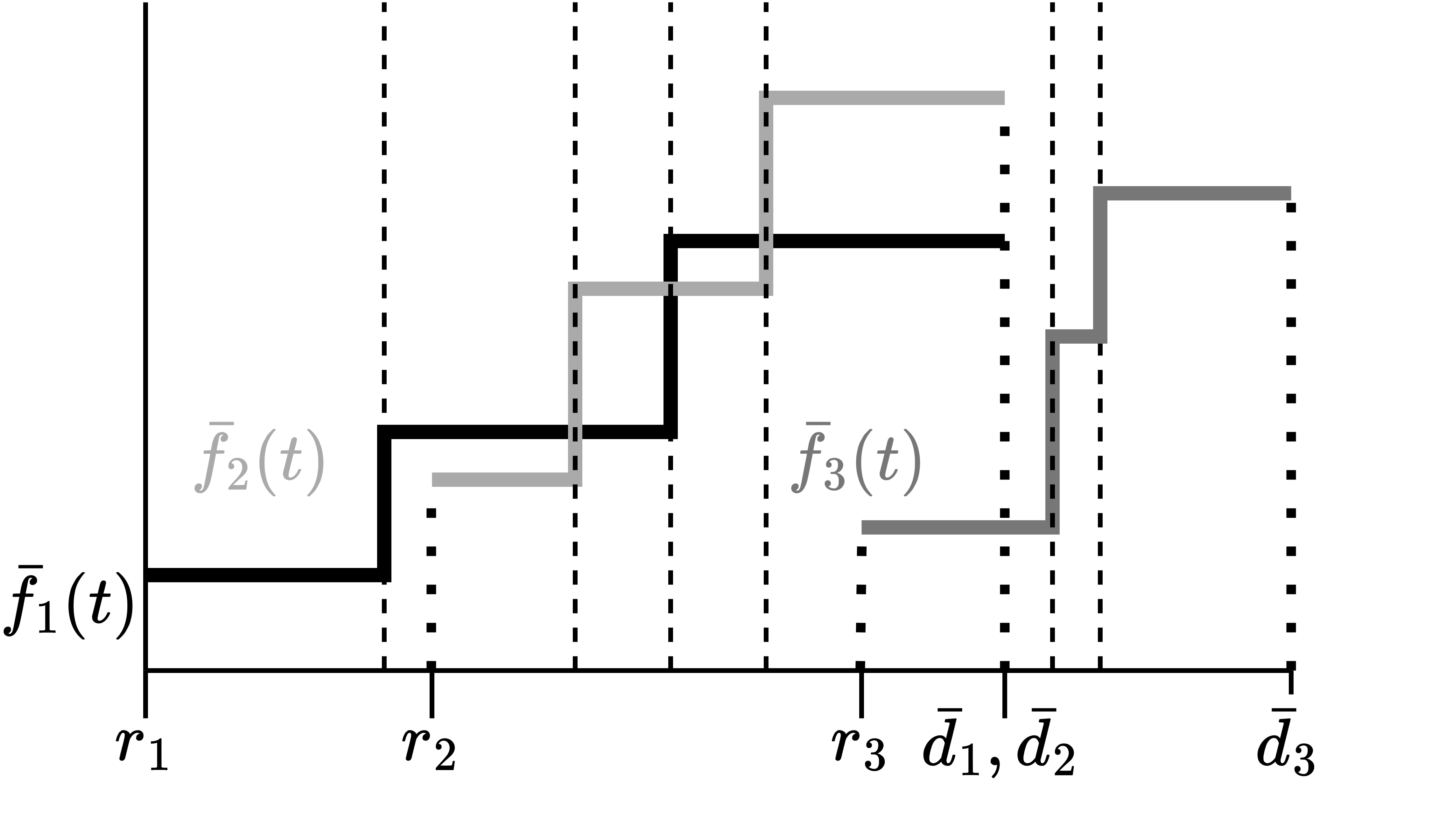}
				\caption{Cost functions}
				\label{fig:gcecsp-example-sub-costfunctions}
			\end{subfigure}%
		}
		\cr
		\noalign{\hfill}
		\hbox{%
			\begin{subfigure}[b]{.45\textwidth}
				\centering
				\includegraphics[width=0.9\textwidth]{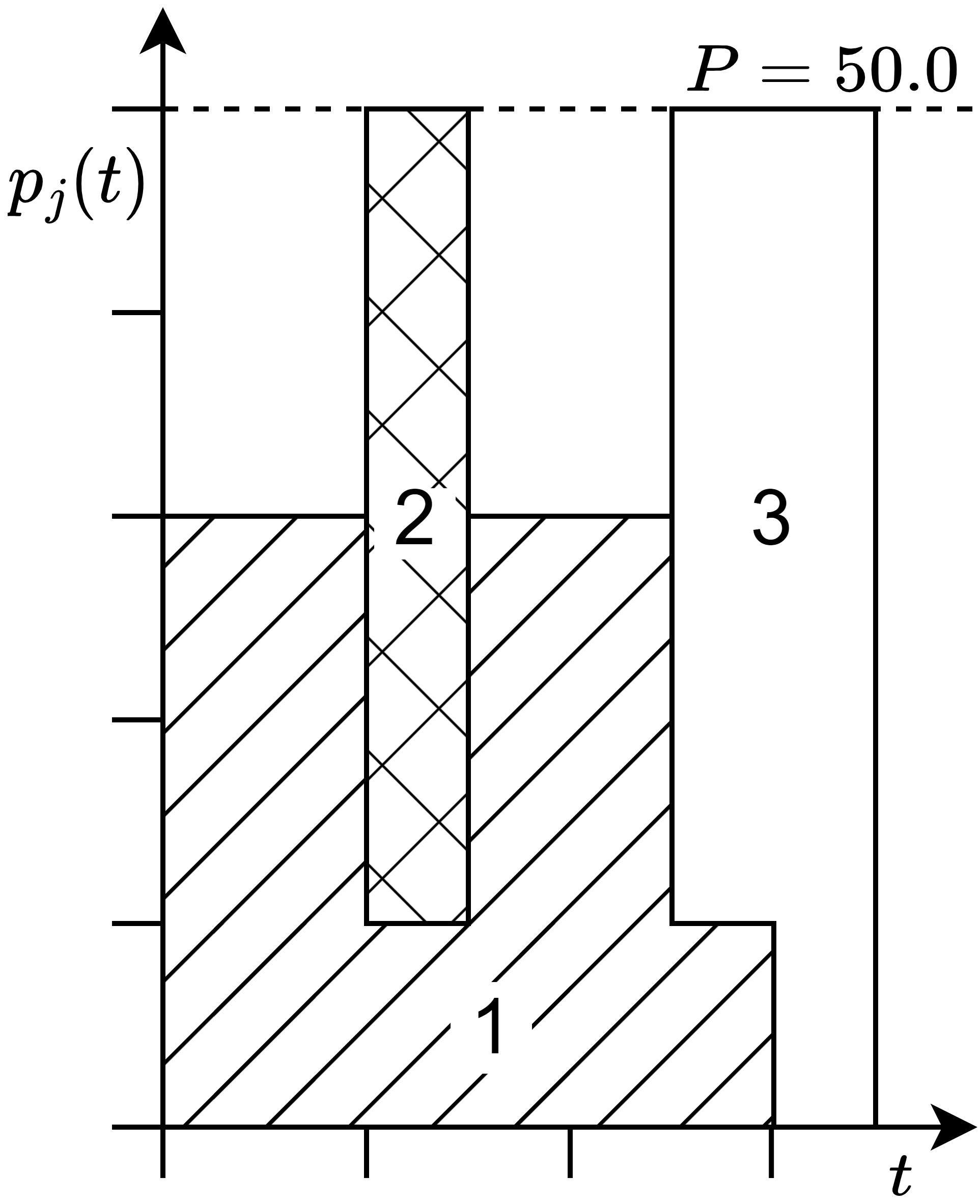}
				\caption{Visualization of a feasible schedule}
				\label{fig:gcecsp-example-sub-schedule}
			\end{subfigure}%
		}\cr
	}
	\caption{Example instance with three jobs}
	\label{fig:gcecsp-example-instance}
\end{figure}

An example instance with three jobs and $P = 50.00$ is provided in Figure \ref{fig:gcecsp-example-instance}. Figure \ref{fig:gcecsp-example-sub-properties} lists the properties of all jobs and Figure \ref{fig:gcecsp-example-sub-schedule} visualizes a feasible schedule for this instance. The cost functions in Figure \ref{fig:gcecsp-example-sub-costfunctions} are step-wise cost functions, which will be discussed in detail in Section \ref{sec:gcecsp-hybrid-alg}.

\section{Hybrid optimization framework} \label{sec:gcecsp-hybrid-alg}

We will introduce our optimization framework by applying it to one particularly interesting variant. In this variant, we want to minimize an increasing step-wise cost function $\bar{f}_j(C_j)$, that determines the cost of completing a job $J_j$ at a given time. Such a function consists of $k$ pieces, each spanning an interval where the cost of completing a job does not change. Essentially, a job now has $k-1$ due dates, or \textit{jump points}, after each of which the cost of completing the job increases, and a single deadline. The jump points are chosen independently for each job. Figure \ref{fig:gcecsp-example-sub-costfunctions} shows example cost functions for the example instance, with two jump points each ($k=3$).

\subsection{Event-based model}\label{sec:gcecsp-sub-model}
Recall that a schedule consists of, for every job $J_j$: a start time $S_j$, a completion time $C_j$ and a resource consumption profile $p_j(t)$. So, for every job, the timing of two \textit{events} has to be determined, which leads to a total of $2n$ events. A sequence of these events $\mathcal{E}$ is the basis of any schedule. We will call these events \textit{plannable}, as we need to decide their exact timing. We denote the (ordered) set of plannable events as $\mathcal{P}$.

\begin{figure}
	\centering
	\includegraphics[width=\textwidth]{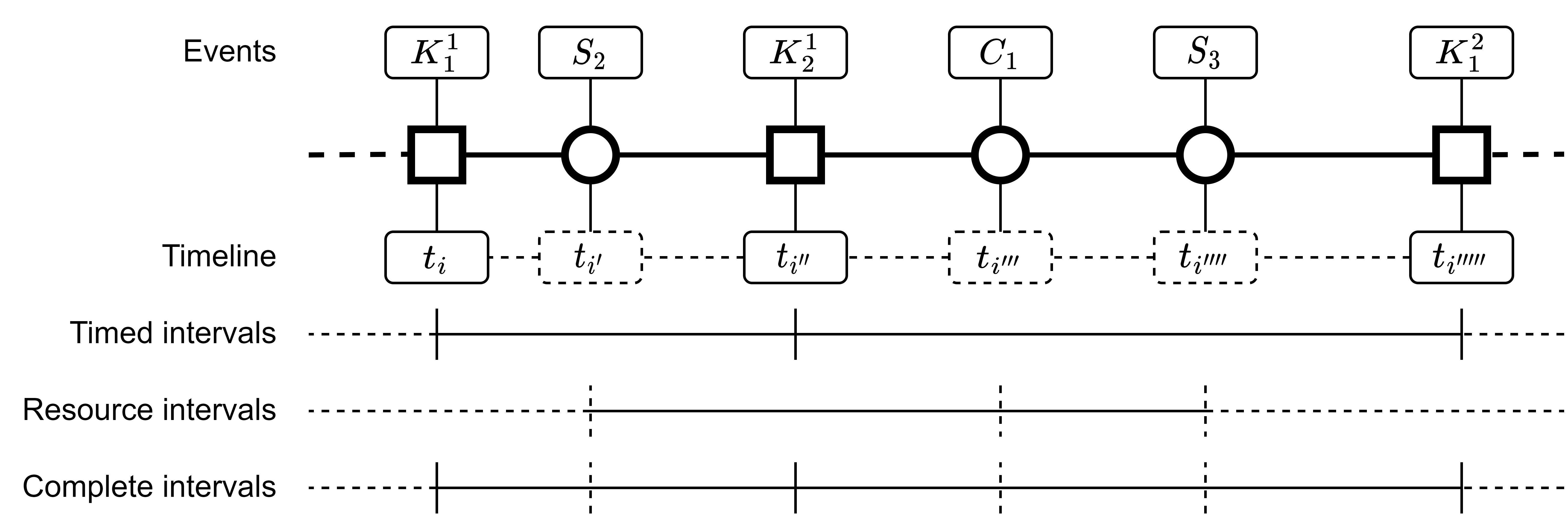}
	\caption{Graphical representation of part of an event order and associated intervals}
	\label{fig:gcecsp-intervals}
\end{figure}

In most cases, we only need these $2n$ events in our model. For the step-wise cost functions, however, we introduce another event type: \textit{fixed-time} events. These are events of which the exact timing is already known. In this case, the $k-1$ jump points $\{K^1_j, ..., K^{k-1}_j\}$ of each job are added to event list $\mathcal{E}$ as fixed-time events, bringing the total to $n(k+1)$ events. We denote the (ordered) set of fixed-time events as $\mathcal{F}$. Including them in the order of events has some benefits for determining the objective value of a potential solution during our algorithm, which will be explained in Section \ref{sec:gcecsp-sub-eval-order}. From a given order of events $\mathcal{E}$, we can identify $2n-1$ \textit{resource intervals}, each defined by two consecutive plannable events in our sequence and $n(k-1) -1$ \textit{timed intervals}, each defined by two consecutive fixed-time events. Part of an example event order and associated intervals is given in Figure \ref{fig:gcecsp-intervals}. Note that, whenever intervals are mentioned without further qualification in the following, these are \textit{resource intervals}. 

The following theorem is a reformulation of the theorem proven by \cite{Brouwer2023} for the case where the objective is to minimize weighted completion time, but the proof can be modified to fit any objective that is a function of (plannable) event time.
\vspace{1em}
\begin{theorem}\label{theorem:gcecsp-constant}
	For any feasible schedule $\mathcal{S}$ with resource consumption profiles $p_j(t)$ for all jobs $J_j, j \in \{1, ..., n\}$ and objective value $W$ that follows a given event order $\mathcal{E}$, a feasible schedule $\mathcal{S}'$ exists with the same order of events $\mathcal{E}$ and objective value $W$ where the resource consumption $p_j'(t)$ of all jobs remains constant during each interval. \qed
\end{theorem}
\vspace{1em}
Theorem \ref{theorem:gcecsp-constant} shows that we can simplify the resource consumption profile function $p_j(t)$ to a discrete number of values $p_{j,i}$, indicating the amount of resource each job $J_j$ consumes during interval $i, i \in \{1, ..., 2n-1\}$. These values $p_{j,i}$ together with event times $t_i$ fully describe the schedule.

The GCECSP can be solved using a mixed-integer linear program, where, besides $t_i$ and $p_{j,i}$, additional binary variables $a_{i,i'}$ and $b_{i,i'}$ are used to model the order of events. We adapt the formulation presented by \cite{Brouwer2023} for the variant with step-wise cost functions. The full MILP formulation is presented below, followed by a description of the model.

\begingroup
\allowdisplaybreaks
\begin{align}
	\begin{aligned}
	\min \sum\limits_{j \in \{1, ..., n\}}w_{j,1}\\
	+ \left( \sum\limits_{l = 2}^{k} w_{j,l} a_{2n + (k-1)(j-1) + (l - 1),2j} \right)
	\end{aligned} & & \text{s.t.} \label{eq:objective}\\
	t_i \leq t_{i'} + Ma_{i',i} & & \forall i, i' \in \{1, ..., (k+1)n\}, i \neq i' \label{eq:event-order}\\
	t_{2j} \leq \bar{d}_j & & \forall j \in \{1, ..., n\} \label{eq:deadline}\\
	t_{2j - 1} \geq r_j & & \forall j \in \{1, ..., n\} \label{eq:release-time}\\
	\sum\limits_{i \in \{1, ..., 2n\}} p_{j,i} = E_j & & \forall j \in \{1, ..., n\} \label{eq:resource-requirement}\\
	\begin{aligned}
	p_{j,i} \geq P_j^- (t_{i'} - t_{i}) - (1-b_{i,i'})M \\
	- (1-a_{2j - 1,i'})M - (1-a_{i, 2j})M
	\end{aligned}
	& & 
	\begin{aligned}
	\forall j \in \{1, ..., n\}, \\
	i, i' \in \{1, ..., 2n\}, i \neq i'
	\end{aligned} \label{eq:lower-bound}\\
	p_{j,i} \leq P_j^+ (t_{i'} - t_{i}) + Ma_{i',i} & &
	\begin{aligned}
		\forall j \in \{1, ..., n\}, \\
		i, i' \in \{1, ..., 2n\}, i \neq i'
	\end{aligned} \label{eq:upper-bound}\\
	p_{j,i} \leq a_{i, 2j}M  & &  \forall j \in \{1, ..., n\}, i \in \{1, ..., 2n\} \label{eq:zero-finish}\\
	p_{j,i} \leq (1- a_{i, 2j - 1})M  & &  \forall j \in \{1, ..., n\}, i \in \{1, ..., 2n\} \label{eq:zero-start}\\
	\sum\limits_{j\in \{1, ..., n\}} p_{j,i} \leq P(t_{i'} - t_{i}) + Ma_{i',i} & & \forall i, i' \in \{1, ..., 2n\}, i \neq i' \label{eq:interval-capacity}\\
	a_{i,i'} + a_{i',i} = 1 & & \forall i, i' \in \{1, ..., (k+1)n\}, i < i' \label{eq:exclusive-order}\\
	\sum\limits_{i'' \in \{1, ..., 2n\}}\left(a_{i,i''} - a_{i',i''}\right) \leq 1 + (1 - b_{i,i'})M & & \forall i, i' \in \{1, ..., 2n\}, i \neq i' \label{eq:one-successor-1}\\
	\sum\limits_{i'' \in \{1, ..., 2n\}}\left(a_{i,i''} - a_{i',i''}\right) \geq 1 - (1 - b_{i,i'})M & & \forall i, i' \in \{1, ..., 2n\}, i \neq i' \label{eq:one-successor-2} \\
	\sum_{i,i'\in \{1, ..., 2n\}, i \neq i'} b_{i,i'} = 2n - 1 & & \label{eq:total-successors}\\
	t_{2j} - t_{2j - 1} \leq E_j / P^-_j & & \forall j \in \{1, ..., n\} \label{eq:max-processing-window}\\
	t_{2j} - t_{2j - 1} \geq E_j / P^+_j & & \forall j \in \{1, ..., n\} \label{eq:min-processing-window}\\
	p_{j,i} \geq 0 & & \forall j \in \{1, ..., n\}, i \in \{1, ..., 2n\} \nonumber \\
	t_{i} \geq 0 & & \forall i \in \{1, ..., 2n\} \nonumber \\
	{\color{gray} t_{i} = c_i} & & {\color{gray}\forall i \in \{2n + 1, ..., (k+1)n\}} \nonumber \\
	a_{i,i'} \in \{0,1\} & & \forall i,i' \in \{1, ..., (k+1)n\}, i \neq i' \nonumber \\
	a_{i,i} = 0 & & \forall i \in \{1, ..., (k+1)n\} \nonumber \\
	b_{i,i'} \in \{0,1\} & & \forall i,i' \in \{1, ..., 2n\}, i \neq i' \nonumber 
\end{align}
\endgroup 

Note that, for the indexing scheme of events, the first $2n$ indices are reserved for plannable events, where $i = 2j - 1$ is the start event of job $J_j$, and $i = 2j$ its completion. The fixed-time events follow, where the $j$th $k-1$ indices (i.e., $2n+(k-1)(j-1)+1, ..., 2n+j(k-1)$) belong to the jump points of job $J_j$. Furthermore, we identify (resource) intervals by the index $i$ of the plannable event that occurs at the start of that interval.

We minimize the sum of the step-wise cost functions (Equation \eqref{eq:objective}), where each jump point that is passed before the completion time adds a given amount to the value of the function. This is subject to the following constraints. Equation \eqref{eq:event-order} controls the timing of events $t_i$, such that it respects the order defined by the order variables $a_{i,i'}$. Equations \eqref{eq:deadline} and \eqref{eq:release-time} enforce the deadline and release time of a job $J_j$. Equation \eqref{eq:resource-requirement} ensures that the resource consumption $p_{j,i}$ of a job $J_j$ over all intervals sums up to its resource requirement $E_j$. Equation \eqref{eq:lower-bound} ensures that the resource consumption $p_{j,i}$ of job $J_j$ in interval $i$ is at least equal to its lower bound $P^-_{j}$ for the duration of the interval if the interval lies between the start and completion events of job $J_j$. Equation \eqref{eq:upper-bound}, in turn, ensures that it does not consume more than its upper bound $P^+_j$. Equations \eqref{eq:zero-finish} and \eqref{eq:zero-start} ensure that the resource consumption $p_{j,i}$ of job $J_j$ is zero outside its processing window. Equation \eqref{eq:interval-capacity} limits summed resource consumption $p_{j,i}$ of all jobs active in interval $i$ to the resource availability $P$ for the duration of that interval. Equation \eqref{eq:exclusive-order} defines the relationship between order variables $a_{i,i'}$ and $a_{i',i}$. Either event $i$ is before $i'$, or $i'$ is before $i$. Equations \eqref{eq:one-successor-1} and \eqref{eq:one-successor-2} together define the successor variables $b_{i,i'}$ in terms of the order variables $a_{i,i'}$. A plannable event $i'$ is a successor of another plannable event $i$ if the difference in the number of plannable events that happen before each event is exactly one. Note that $b_{i,i'}$ is defined on plannable events only, as the successor relationship only matters for the enforcement of lower bounds in resource intervals. Equation \eqref{eq:total-successors} states that there are exactly $2n-1$ successor relationships among the $2n$ plannable events. Finally, Equations \eqref{eq:max-processing-window} and \eqref{eq:min-processing-window} are valid inequalities that model a maximum and minimum duration for the processing window of a job $J_j$.

\subsection{Decomposition and local search} \label{sec:gcecsp-sub-ls}
We observe that, after fixing the binary variables in the MILP, we are left with an LP that quickly finds an optimal schedule for a given order of events. This inspires the decomposition of the problem in two parts: (1) finding a good event order $\mathcal{E}$, and (2) determining an optimal schedule for a given event order. This idea forms the basis of the approach that we developed in our previous work \parencite{Brouwer2023}, where we use simulated annealing to explore permutations of the order of events $\mathcal{E}$, with neighborhood operators that swap two adjacent events and occasionally move an event multiple positions forwards or backwards in the event order. In every iteration of the simulated annealing, we evaluate the objective by solving an LP to find the optimal schedule for the given event order, defining event times $t_i$ and resource consumption profiles $p_{j,i}$.

This approach generalizes well to the GCECSP. It is a shared property of the GCECSP and its variants that it is difficult to find a good event order $\mathcal{E}$, while finding a schedule for a given event order is relatively easy. The general approach can be applied to the GCECSP and its variants, while some elements, such as the neighborhood operators or the means of assessing the quality of event orders, may need to be tailored in each case. This tailoring, however, may also open up opportunities for increasing the performance of the approach. This is the case for the variant with step-wise cost functions, as will become clear in the remainder of this work.

In the remainder of this section we will discuss how to tailor the general approach for the variant with step-wise cost functions. We need to make a number of changes: (1) we adapt the greedy algorithm for finding an initial solution; (2) we extend the set of implicit precedence constraints that we consider in the local search; (3) because of the increased number of events for an instance of the same size, when compared to the general case, we adapt the neighborhood operator to generate neighbors by moving an event multiple positions more often; (4) to avoid spending a lot of time on a plateau, we terminate the search if the number of iterations since the last improvement of the current solution exceeds one cooling period (see Section \ref{sec:gcecsp-sub-instances}); (5) we change the way the quality of event orders is assessed, which will be explained in Section \ref{sec:gcecsp-sub-eval-order}. Because we use fixed-time events to model jump points in the cost function, the value of the objective function is fully determined by the event order. However, we still need to determine whether a feasible solution exists for that particular event order.

\subsection{Assessing the quality of event orders}\label{sec:gcecsp-sub-eval-order}
In the general approach, an LP is used to solve the subproblem, i.e. for determining an optimal schedule for the given event order (if one exists). In order to give the simulated annealing search enough flexibility, we allow it to search through the space of infeasible solutions. To evaluate infeasible candidate solutions, we include three types of so-called \textit{violation variables} ($s^+_{j,i}, s^-_{j,i}, s^t_{i} \geq 0$) that allow for the violation of the constraints enforcing the (1) upper and (2) lower bound on the resource consumption for a job, and (3) the upper bound on the amount of available resource during an interval at the cost of a penalty in the objective. If the sum of the violation variables is larger than 0, no feasible solution exists for that order, but we can still view a decrease in the penalty term as a move in the right direction.

While the search is ongoing, however, it is sufficient if we know the objective value, including the penalty term. We do not need the schedule, i.e. the exact values of all event times $t_i$ and resource consumption during intervals $p_{j,i}$. For the case with step-wise cost functions, this means that we may not have to solve the LP every time. The value of the objective (excluding the penalty term) is fully determined by the order of events. For each completion event $C_j$ in the event order, let $K_j^{\rightarrow}$ denote the next jump point of job $J_j$ in the event order (where $K_j^{\rightarrow} = \bar{d}_j$ if no further jump point occurs). Then, any \textit{feasible} schedule with the given event order will have a cost of exactly $\sum_j \bar{f}_j(K_j^{\rightarrow})$.

At this point, we exploit that this base score may already be sufficient to reject a candidate solution. We pre-compute the maximum objective value that we will accept during the evaluation. If the base score already exceeds this, we will not compute the penalty term at all.

For estimating the penalty term, we propose three alternatives:
\begin{enumerate}
	\item solve the LP, which will always yield the true value of the penalty term;
	\item compute a maximum flow on a bipartite graph, modeling a relaxed version of the problem;
	\item compute lower bounds on the value of the violation variables using only the event order.
\end{enumerate}

The results of extensive tests of using each of the three methods for computing the penalty term are presented in Section \ref{sec:gcecsp-sub-instances}.

\subsubsection{Solving the LP}

The LP is, as before, what remains if we fix the order variables in the MILP and add violation variables to the relevant constraints. The objective is then to minimize the weighted sum of the violation variables.

We show the full LP formulation below. For ease of readability, we introduce some additional notation:
\begin{itemize}
		\item $X_{j, i} = \begin{cases}
		1 & \text{if } I(t_{2j-1}) \leq I(t_i) \text{ and } I(t_{2j}) \geq I(t_i) \\
		0 & \text{otherwise}
	\end{cases}$ defines if a job $J_j$ is active during interval $i$, which we already know, based on the order of events;
	\item $\text{succ}(i) = E(I(i) + 1)$, the index of the successor of event $i$ in the order of events $\mathcal{E}$;
	\item $t_{\text{succ}(i)} = \infty$ for the last event in $\mathcal{E}$;
	\item $\text{succ}_{\mathcal{P}}(i)$ is the same as $\text{succ}(i)$, but skips over fixed-time events;
	\item $p_{j,i} = 0, \forall j$ for the last event in the order of events $\mathcal{E}$.
\end{itemize}
\begingroup\allowdisplaybreaks
\begin{align*}
	\min \sum\limits_{i = 1}^{2n} \left( L^R s^t_i + \sum\limits_{j = 1}^n L^B(s^-_{j,i} + s^+_{j,i}) \right) & & \text{s.t.} \\
	t_i \leq t_{\text{succ}(i)} & & \forall i \in \{1, ..., (k+1)n\} \\ 
	t_{2j} \leq \bar{d}_j & & \forall j \in \{1, ..., n\} \\ 
	t_{2j - 1} \geq r_j & & \forall j \in \{1, ..., n\}  \\ 
	\sum\limits_{i = 1}^{2n} p_{j,i} = E_j & & \forall j \in \{1, ..., n\} \\ 
	p_{j,i} \geq P_j^- X_{j,i} (t_{\text{succ}_{\mathcal{P}}(i)} - t_{i}) - s^-_{j,i} & &  \forall j \in \{1, ..., n\}, i \in \{1, ..., 2n\} \\ 
	p_{j,i} \leq P_j^+ X_{j,i} (t_{\text{succ}_{\mathcal{P}}(i)} - t_{i}) + s^+_{j,i} & &  \forall j \in \{1, ..., n\}, i \in \{1, ..., 2n\} \\ 
	\sum\limits_{j = 1}^n p_{j,i} \leq P(t_{\text{succ}_{\mathcal{P}}(i)} - t_{i}) + s^t_i & & \forall i \in \{1, ..., 2n\} \\ 
	t_{i} \geq 0 & & \forall i \in \{1, ..., 2n\} \nonumber \\
	t_{i} = c_i & & \forall i \in \{2n + 1, ..., (k+1)n\} \nonumber \\
	p_{j,i} \geq 0 & & \forall j \in \{1, ..., n\}, i \in \{1, ..., 2n\} \nonumber \\
	s^-_{j,i} \geq 0 & & \forall j \in \{1, ..., n\}, i \in \{1, ..., 2n\} \nonumber \\
	s^+_{j,i} \geq 0 & & \forall j \in \{1, ..., n\}, i \in \{1, ..., 2n\} \nonumber \\
	s^t_{i} \geq 0 & & \forall i \in \{1, ..., 2n\} \nonumber
\end{align*}
\endgroup

\subsubsection{Max flow on a bipartite graph}
We can model the assignment of resources to jobs and intervals as a flow problem on a bipartite graph, if we relax some constraints. The bipartite graph models the given event order by two layers of nodes. On one side, every node represents a job, and on the other, every node represents the interval between two consecutive fixed-time events. An arc between the two exists if the job is active during that interval, i.e. if the start event of the job happens in or before that interval and its completion event happens in or after it. We set the arc capacity such that upper and lower bounds are enforced wherever possible. This will be explained in more detail below.

Arcs from the source to a job node representing job $J_j$ initially have capacity equal to the corresponding resource requirement $E_j$, arcs from a job node for job $J_j$ to an interval node $[\mathcal{F}_i,\mathcal{F}_{i'}]$ have capacity $(\mathcal{F}_{i'} - \mathcal{F}_{i})P^+_j$, enforcing the upper bound, and arcs from an interval node $[\mathcal{F}_i,\mathcal{F}_{i'}]$ to the sink have capacity $(\mathcal{F}_{i'} - \mathcal{F}_{i})P$, modeling the available resource in that interval. Lower bounds can be enforced on intervals that do not contain the start or completion event of a job, resulting in the adjustment of arc capacities on the arcs on the path from the source to the sink going through the specific job and interval node. The capacity of each of the arcs on this path is reduced by the amount of resource required to guarantee the lower bound for the job in that interval (i.e. $(\mathcal{F}_{i'} - \mathcal{F}_{i})P^-_j$).

\begin{figure}
	\centering
	\includegraphics[width=\textwidth]{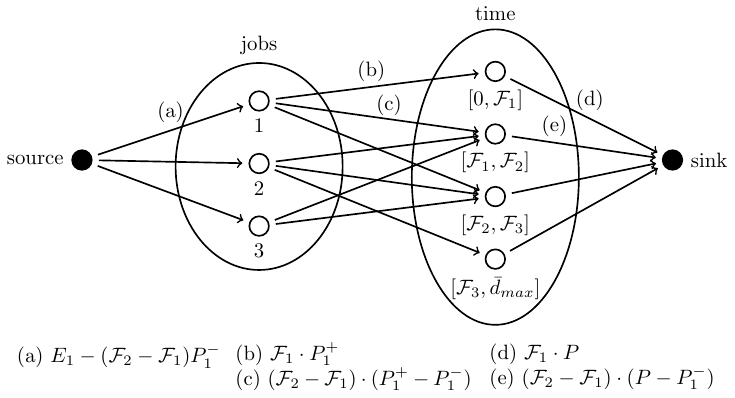}
	\caption{Example with three jobs and four intervals}
	\label{fig:gcecsp-flow-check}
\end{figure}

An example graph is shown in Figure \ref{fig:gcecsp-flow-check}. The capacities of a few arcs have been added in the figure. In this example, job 1 is active during the first three intervals. Therefore, we know it is active during the entire second interval, and we can already assign the lower bound. This results in the adjustment of the capacity on the arcs that are on the (unique) path from source to sink that goes through the nodes representing job 1 and interval 2 (a, c and e).

The assignment of lower bounds may already lead to the conclusion that the event order is infeasible. If it results in negative capacities in the arcs from source to job nodes, we know that a job violates its lower bounds by at least that amount. We keep track of this, and put the capacity to 0. This will contribute to the estimation of the penalty function, analogous to the $s^-_{j,i}$ variable in the LP.

Similarly, if an arc from a time node to the sink has a negative capacity, we know that there is insufficient resource availability to let all jobs consume at least their lower bound. Again, we keep track of this as a minimum violation of the resource availability, analogous to the $s^t_{i}$ variable in the LP.

We compute the max flow through the resulting network. The difference between the value of the max flow and the sum of the (adjusted) resource requirements gives an indication of the total resource shortage, and therefore of the value of the resource violation ($\sum_i s^t_i$). We add to this the penalties that we already encountered while computing the arc capacities.

If the difference between these two values is zero, it is a strong indication that a feasible solution exists, or an infeasible solution (due to the order of start and completion events in a fixed-time interval, in combination with restrictive lower or upper bounds) that is fairly close to being feasible. The smaller the expected number of plannable events within a fixed-time interval, the less often such a flow will result from an infeasible event order.

A feasible flow does not guarantee a feasible solution. It does not enforce the order of plannable events within an interval in between two fixed-time events. Therefore, the resulting flow may not be valid for the exact event order provided. Ensuring one plannable event to respect its bounds may, due to the order of events, push another to a time within the interval where it cannot respect its bounds any more. The estimation that results from this flow is almost always a close underestimation. However, it is not a true lower bound on the penalty term, as it does not (properly) account for the possible violation of lower and upper bounds ($s^-_{j,i}$ and $s^+_{j,i}$ in the LP).

\subsubsection{Simple lower bounds}
Lower bounds can be computed on each of the three types of violation variable. For the purpose of these lower bounds, we will introduce $\mathcal{F}^{\leftarrow}(i)$ and $\mathcal{F}^{\rightarrow}(i)$, indicating the closest fixed-time event preceding or following an event $i$, respectively. These functions can be constructed by going over the order of events $\mathcal{E}$ once in both directions, and efficiently updated every time an event is moved.

We can determine a lower bound on the violation of the lower bounds in the following way:
A job $J_j$ can start no later than the first fixed-time event that follows its start event in the event order: $\mathcal{F}^{\rightarrow}(S_j)$, and complete no sooner than the last fixed-time event that preceded its completion event: $\mathcal{F}^{\leftarrow}(C_j)$. The time between these two time points is the minimum processing time. We can use this to determine the minimum amount of resource that the job consumes if it respects its lower bound during this time. If this amount is larger than its total resource requirement, the difference between the two gives us a lower bound on the value of $\sum_i s^-_{j,i}$:

\[
\sum_{j=1}^n\max\left(0, \left(t_{\mathcal{F}^{\leftarrow}(C_j)} - t_{\mathcal{F}^{\rightarrow}(S_j)}\right)P^-_j - E_j\right)
\]

Similarly, a lower bound on the violation of the upper bound for a job $J_j$ ($s^+_{j,i}$) can be computed using the maximum processing time, based on the last fixed time event preceding the start event $\mathcal{F}^{\leftarrow}(S_j)$ and the first fixed time event following the completion event $\mathcal{F}^{\rightarrow}(C_j)$:

\[
\sum_{j=1}^n\max\left(0, E_j - \left(t_{\mathcal{F}^{\rightarrow}(C_j)} - t_{\mathcal{F}^{\leftarrow}(S_j)}\right)P^+_j\right)
\]

Finally, it remains to estimate the value of the violation variable controlling the amount of available resource during an interval. We go over the event list once, adding the resource requirement $E_j$ to the total of consumed resource up to that point if we encounter a completion event $C_j$. With every fixed-time event, we compute the total amount of resource available up to that point ($P\cdot t_i$, where $i$ is the index of the fixed-time event). Any time the newly computed amount of available resource is lower than the total amount of consumed resource, we add the difference to our estimation of the violation. At the end of the event list, this results in a lower bound on the resource shortage during intervals.
\begin{algorithm}[H]
\caption{Estimate interval capacity penalty}
\begin{algorithmic}[1]
	\State $E_{total} \gets 0$
	\State $E_{shortage} \gets 0$
	\For{$e \in \mathcal{E}$}
		\If{IsCompletionEvent($e$)}
			\State $E_{total} \gets E_{total} + E_{e/2}$
		\EndIf
		\If{IsFixedTimeEvent($e$)}
			\State $E_{shortage} \gets E_{shortage} + \max(0, E_{total} - \text{t}_{I(e)} \cdot P)$
			\State $E_{total} = \min{E_{total}, \text{t}_{I(e)} \cdot P}$
		\EndIf
	\EndFor
	\State \Return $E_{shortage} \cdot L^R$\;
\end{algorithmic}
\end{algorithm}

The weighted combination of these three lower bounds provides a lower bound on the total value of the penalty term. All three can be computed in $O(n)$ time. Storing intermediate results allows for efficient updates in an amount of time linear in the number of positions that an event is moved.

\section{Applicability of the framework for the general problem} \label{sec:gcecsp-variants}
In this section, we aim to give an overview of the applicability of our approach: we will discuss the adaptations needed to apply our framework to an elaborate set of variants of the GCECSP. Our main aim is to show the flexibility and adaptability of our framework, where we do not claim that our approach is the best choice in all cases where it can be applied.

We will go through elements of the problem in three sections, to discuss possible variants and extensions, first discussing the range of possible objective functions (Section \ref{sec:gcecsp-sub-objectives}), then properties of jobs (Section \ref{sec:gcecsp-sub-job-props}) and finally the resource (Section \ref{sec:gcecsp-sub-resources}).

\subsection{Objectives}\label{sec:gcecsp-sub-objectives}
In general, the approach can easily be adapted to work for any objective that is a linear function of variables that are present in the LP. For objectives that are a function of completion time $C_j$ (or more general, any plannable event time $t_i$) we distinguish between three types that we can model with increasing effort. We will discuss these three types, with a number of examples, followed by some examples that do not belong to these three types. Then, we will briefly discuss objectives that are not a function of event times.

\subsubsection{Type I: Linear functions}
For an objective that is a linear function of event times, only the replacement of the objective in the LP is required. Examples of this type of objectives include:
\begin{itemize}
	\item Minimize (weighted) completion time $\sum_j w_jC_j$ (see our previous work \parencite{Brouwer2023});
	\item Minimize total processing time $\sum_j (C_j - S_j)$.
\end{itemize}

\subsubsection{Type II: Minmax criteria (and extensions)}
These objectives require a minimal amount of additional adaptations to the LP compared to objectives of Type I. In addition to the change of the objective, for example, a helper variable and some additional constraints may have to be added to compute the objective. Examples of these are objectives that contain a min or max function, such as minimizing the makespan. They can be modeled using a helper variable for the objective ($\min w$) and a number of linear equations ($\forall j: w \geq C_j$), but not in a single linear expression. Examples of this type of objectives include:
\begin{itemize}
	\item Minimize makespan $\max_j\{C_j\}$;
	\item Minimize maximum lateness $\max_j\{C_j - d_j\}$;
\end{itemize}

\subsubsection{Type III: Piece-wise linear functions} Objectives that are piece-wise linear can be modeled by splitting them in their constituent pieces using additional fixed-time events, as discussed in Section \ref{sec:gcecsp-sub-model}. For this type of objectives it is not sufficient to adapt the LP (as it is for Types I and II), but changes to the event model, or the introduction of binary variables in our LP are required. Examples of this type of objectives include:
\begin{itemize}
	\item Minimize (weighted) number of tardy jobs $\sum_j w_jU_j$, $U_j = \begin{cases} 0 & C_j \leq d_j \\ 1 & \text{otherwise}\end{cases}$;
	\item Minimize a step-wise function of completion time $\sum \bar{f_j}(C_j)$ (see Section \ref{sec:gcecsp-hybrid-alg}).
\end{itemize}

\subsubsection{Non-linear functions}
If the objective is a non-linear function of event times, the LP-solver can be replaced by an NLP approach.

\subsubsection{Functions of resource variables} Objectives that are a function of resource consumption $p_{j,i}$, can be accommodated along the same lines as functions of event times, as long as Theorem \ref{theorem:gcecsp-constant} holds for the particular objective. This is the case if a solution with constant resource consumption during intervals dominates solutions with other consumption profiles. Some examples of such functions are:

\begin{itemize}
	\item Type I: minimize total resource consumption: $\sum_{j,i} p_{j,i}$, where jobs have \textit{efficiency functions} (see Section \ref{sec:gcecsp-sub-job-props});
	\item Type III: minimize resource cost $\sum_{j,i} c_{i,j} p_{j,i}$, where cost fluctuates across predetermined intervals. Note that the complete set of intervals have to be used in the model, where an interval is defined as the time between two consecutive events of \emph{any} type, instead of just resource intervals (see Section \ref{sec:gcecsp-sub-model});
	\item Non-linear: any objective using the consumption rate $p_{j,i} / (t_{i'} - t_i)$.
\end{itemize}

\subsubsection{Order-determined objectives}
For any objective that can already be computed based on the order of events alone, the LP is only used to determine the feasibility of the given order of events. Any (efficient) computation of the objective can be performed in the evaluation of a candidate solution before solving the LP. An example of such an objective is the minimization of the maximum number of jobs that are being processed in parallel.

\subsection{Job properties}\label{sec:gcecsp-sub-job-props}
In the description of the GCECSP in Section \ref{sec:gcecsp-prob-desc} we listed a number of properties for each job. We will discuss the impact on our approach of the absence of some of these properties, and the addition of others. The removal of \textbf{release times}, \textbf{deadlines} and \textbf{upper bounds} on the resource consumption rate are fairly straight-forward, as is the introduction of jobs with a \textbf{fixed processing rate}. We will not discuss these adaptations in detail.

The \textbf{lower bound} on the resource consumption rate is an essential property of this type of problem. Maintaining the flexibility of the jobs while removing the lower bound introduces a form of preemption. The problem becomes significantly easier. While the presented approach can still be used, using a flow-based or greedy algorithm will find the optimal solution faster.

In general, \textbf{preemption} can only be modeled to a limited extent. Removing the lower bound allows preemption, but if a lower bound needs to be maintained, allowing preemption can only be done by modeling a single task as multiple jobs with precedence constraints. Then, preemption is allowed a fixed number of times. General preemption can be approximated by making this number large, but this hurts the efficiency of the approach.

\textbf{Efficiency functions} that apply to the conversion of the amount of consumed resource to the contribution towards the resource requirement of a job can be easily integrated in the LP, as long as these functions are linear.

The current approach already considers implicit \textbf{precedence} relations, that follow from the release times and deadlines of jobs. It can easily be extended to include explicit precedence relations between any two events. The local search will not explore any candidate solution that does not respect such a relation.

\textbf{Synchronization} constraints can be included by changing the interpretation of events. In the event order we treat each pair of events that has to happen simultaneously as a single event. Moreover, we add an equality constraint for their event times. In this way, it is ensured that any considered solution will synchronize these two events.

\subsection{Resource}\label{sec:gcecsp-sub-resources}
In the GCECSP, we have a single constant renewable resource. Variations on this aspect of the problem can also be modeled within our framework. We will discuss two example extensions.

Any piece-wise linear \textbf{resource availability function} can be modeled using fixed-time events similar to piece-wise linear objective functions (see Section \ref{sec:gcecsp-sub-objectives}).

\textbf{Multiple resources} can be modeled using additional resource variables $p'_{j,i}$. The bounds, resource requirement, and events can be unique to each resource, or shared among multiple resources. 

A \textbf{discrete resource} can be approximated by using multiple (synchronized) jobs with identical lower and upper bounds to model a single task. This will, however, limit the number of times the resource consumption profile can change. In addition, it will hurt the efficiency of the approach.

\section{Computational results} \label{sec:gcecsp-results}
We will show results for the GCECSP with step-wise cost functions (see Section \ref{sec:gcecsp-hybrid-alg}), using three variants of our approach. For small instances, we will compare these approaches to the MILP (solved with a time limit of 3600 seconds) as well. The data and code used are available online \parencite{DataBrouwer2023, SoftwareBrouwer2023}.
In Section \ref{sec:gcecsp-sub-instances}, we discuss the instances and parameter settings used in our evaluation. Following that, we present the variants of our approach used in the tests and evaluate their performance in Section \ref{sec:gcecsp-sub-results}.
The full result tables are provided in Appendix \ref{app:gcecsp-results}.

\subsection{Test instances and parameter settings}\label{sec:gcecsp-sub-instances}
For generating instances, we will use the procedure developed for instances for  the CECSP in our previous work \parencite{Brouwer2023}. Only the generation of weights will be different, and the generation of times for the jump points has to be added. Recall that $k$ indicates the number of cost intervals in the cost function for each job.

Instead of a single weight $w_j$ from $U(0, 5)$ we generate an array of $k$ weights, in the following way: we first generate an initial weight $w_{j,1}$ from $U(0, 5)$ and then generate $k-1$ weights from $U(w_{j,1}, 5)$. The resulting array of $k$ weights is then sorted to be increasing. Finally, the values are modified to be incremental, such that the cost of completing job $J_j$ right after the first jump point becomes $w_{j,1} + w_{j,2}$, the cost of completing the job right after the second jump point becomes $w_{j,1} + w_{j,2} + w_{j,3}$, and so on.

Then, we generate $k-1$ jump points for each job. The release time $r_j$ and deadline $\bar{d}_j$ are generated as before. We require $k-1$ additional time points in between where the weight changes. We sample these $k-1$ values from $U(r_j + \frac{E_j}{P^+_j}, \bar{d}_j)$ and put them in increasing order.

This completes the modified instance generation approach. Using this, we generate a set of 72 instances: 4 unique instances for each combination of $n$ and $k$ for $n \in \{5, 10, 15, 20, 30, 50\}, k \in \{2, 3, 4\}$.

\begin{figure}
	\centering
	\includegraphics[width=\textwidth]{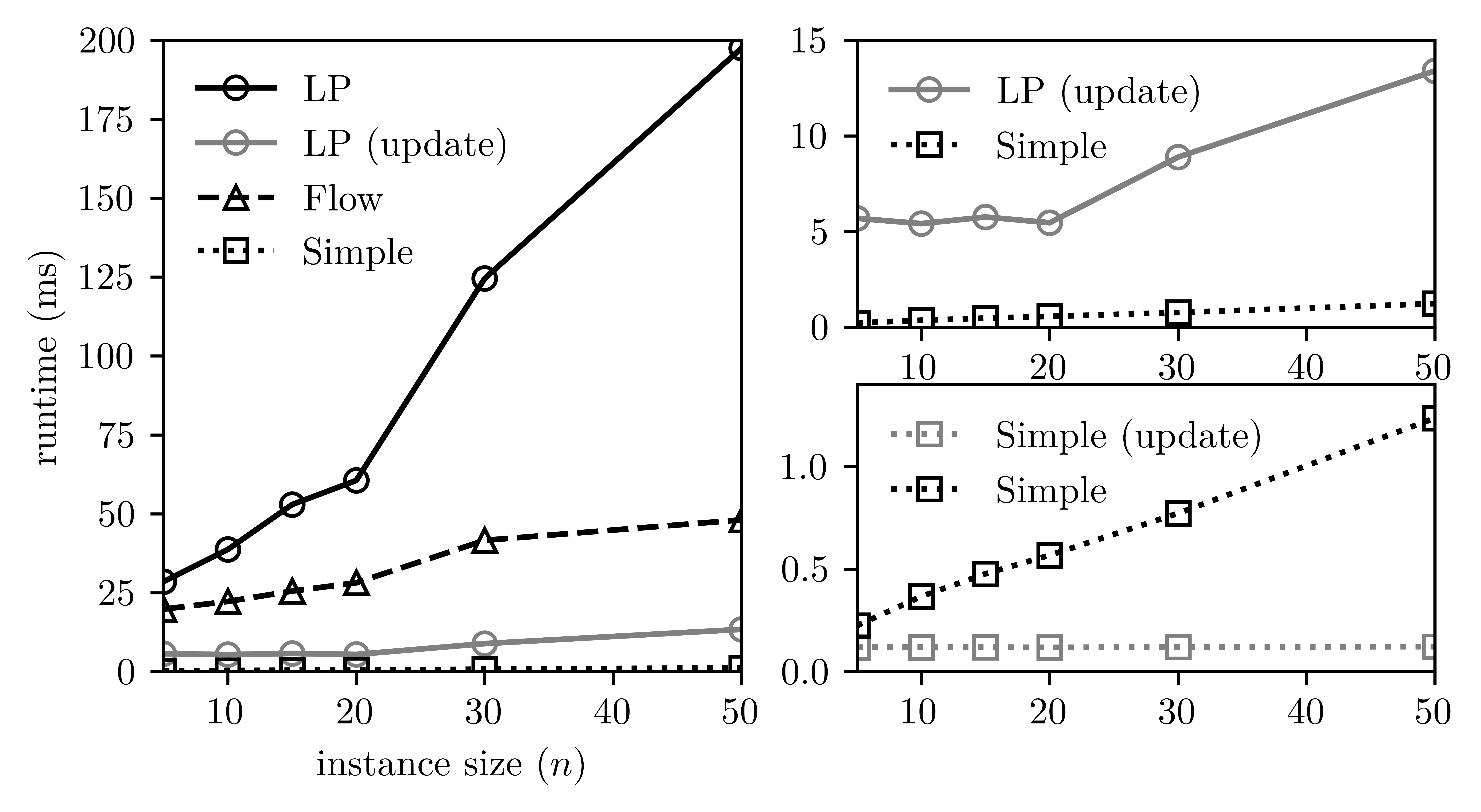}
	\caption{Average runtimes of evaluation methods for event orders} \label{fig:gcecsp-compare-bounds}
\end{figure}

We evaluated the speed and quality of the estimation methods described in Section \ref{sec:gcecsp-sub-eval-order}. For all instance sizes, we ran a test where we computed the estimation methods on a sequence of event orders. The sequence was obtained by performing a random walk of length 1000, randomly shifting an event by one position each time. Averages of the run times are displayed in Figure \ref{fig:gcecsp-compare-bounds}.

The lower bound estimation is much faster than either of the other methods and displays linear scaling behavior. In terms of quality, we note that the computed bounds are not tight. The flow typically results in a better estimation of the penalty term. However, for guiding the search, a strong correlation with the actual penalty term is more important than a tight estimate. To evaluate this, we computed Spearman's rank correlation of the lower bounds and the penalty term, for the instance sizes where $n\geq10$. These all have values between 0.67 and 0.86 (0.77 on average), showing a strong correlation. The flow estimation shows an even stronger correlation (0.82 on average). However, we decided to drop it from our approach, given the limited improvement on the quality of the estimation and the increased runtime compared to the `simple' bounds.

Each of these estimates can be updated efficiently from a previous event order. We have implemented this for both the LP (1) and the `simple' lower bound estimation (3). These are used in the approaches that will be evaluated in Section \ref{sec:gcecsp-sub-results}.

To make sure that we get the best performance out of the presented approaches, we carefully select values for a number of parameters. We chose a commonly used value for the multiplication factor for updating the temperature $\alpha = 0.95$. We performed a limited grid search using different values for the initial temperature $T_{\mathrm{init}} \in \{0.1, 0.2, 0.5, 1, 2, 5\}n$, the number of iterations between temperature updates $\alpha_{\mathrm{period}} \in \{2, 4, 6, 8\}2n$ and the multipliers used the penalize violation variables $L^R, L^B \in \{0.5, 1, 2, 5\}$. We performed runs for all combinations on twelve instances with ten jobs. We observe that finding a good initial solution is important, and exploring too much into a part of the search space consisting of infeasible solutions is detrimental to the eventual success. Our results indicate that it is beneficial to select an initial temperature that is not too high, while cooling down more slowly, and a large penalty for the violation variables, directing the search towards feasible solutions. The addition of a tabu list of length 1, forbidding the movement of the same event two times in a row, yields a slight improvement. The final parameter settings are as follows: $T_{\mathrm{init}} = 0.2n$; $\alpha = 0.95$; $\alpha_{\mathrm{period}} = 16n$; violation penalty: $5.0$.

\subsection{Results and discussion}\label{sec:gcecsp-sub-results}
\begin{table}
	\centering
	\begin{minipage}[b]{0.5\linewidth}
	\begin{tabular}{@{}l|rrr@{}}
		\toprule
		\textbf{Approach} & \textbf{feas.} & \textbf{best} & \textbf{dist.} \\
		\midrule
		\multicolumn{4}{l}{$n=5$} \\
		\midrule
		MILP & 12/12 & 12/12 & 0.00 \\
		SA-LP & 12/12 & 5/12 & 3.48 \\
		SA-MIX & 11/12 & 8/12 & 1.05 \\
		SA-2PH. & 12/12 & 7/12 & 4.45 \\
		\midrule
		\multicolumn{4}{l}{$n=10$} \\
		\midrule
		MILP & 12/12 & 9/12 & 1.06 \\
		SA-LP & 12/12 & 5/12 & 1.88 \\
		SA-MIX & 10/12 & 1/12 & 7.23 \\
		SA-2PH. & 11/12 & 3/12 & 2.76 \\
		\midrule
		\multicolumn{4}{l}{$n=15$} \\
		\midrule
		MILP & 7/12 & 1/12 & 4.98 \\
		SA-LP & 12/12 & 8/12 & 0.38 \\
		SA-MIX & 11/12 & 0/12 & 19.94 \\
		SA-2PH. & 12/12 & 6/12 & 2.18 \\
		\bottomrule
	\end{tabular}
	\end{minipage}
	\quad
	\begin{minipage}[b]{0.5\linewidth}
	\begin{tabular}{@{}l|rrr@{}}
		\toprule
		\textbf{Approach} & \textbf{feas.} & \textbf{best} & \textbf{dist.} \\
		\midrule
		\multicolumn{4}{l}{$n=20$} \\
		\midrule
		SA-LP & 12/12 & 4/12 & 0.98 \\
		SA-MIX & 4/12 & 0/12 & 11.02 \\
		SA-2PH. & 12/12 & 8/12 & 0.21 \\
		\midrule
		\multicolumn{4}{l}{$n=30$} \\
		\midrule
		SA-LP & 12/12 & 5/12 & 0.89 \\
		SA-MIX & 10/12 & 0/12 & 14.06 \\
		SA-2PH. & 12/12 & 7/12 & 1.33 \\
		\midrule
		\multicolumn{4}{l}{$n=50$} \\
		\midrule
		SA-LP & 12/12 & 7/12 & 0.69 \\
		SA-MIX & 2/12 & 0/12 & 10.41 \\
		SA-2PH. & 12/12 & 5/12 & 1.04 \\
		\bottomrule
		\multicolumn{4}{l}{ } \\
		\multicolumn{4}{l}{ } \\
		\multicolumn{4}{l}{ } \\
	\end{tabular}
	\end{minipage}
	\caption{Summary of results for the presented approaches}\label{tab:gcecsp-results}
\end{table}
We present results for the following three variants of our hybrid approach: 
\begin{description}
	\item[SA-LP] uses only the (iteratively updated) LP for the subproblem;
	\item[SA-MIX] uses the `simple' bound estimations to assess the quality of an event order, but also computes the LP if the result of the estimation is 0;
	\item[SA-2PHASE] uses the `simple' bound estimations until it finds a solution for which the estimation is 0, and switches to the LP afterwards.
\end{description}

As discussed in Section \ref{sec:gcecsp-sub-instances}, we will use the `simple' lower bound estimation and the LP to assess the quality of candidate solutions. Preliminary tests showed that using only the `simple' lower bound estimation does not result in good solutions. Therefore, we evaluate approaches combining the use of the `simple' lower bound estimation with the LP.

In Table \ref{tab:gcecsp-results}, we summarize the results for each instance size. For every approach, we list three measures of its performance: the number of instances for which it found a feasible solution, the number of instances for which it found the best solution among the tested approaches, and the average distance (in \%) to the best solution among the tested approaches (infeasible solutions are excluded).

\begin{figure}
	\centering
	\includegraphics[width=\textwidth]{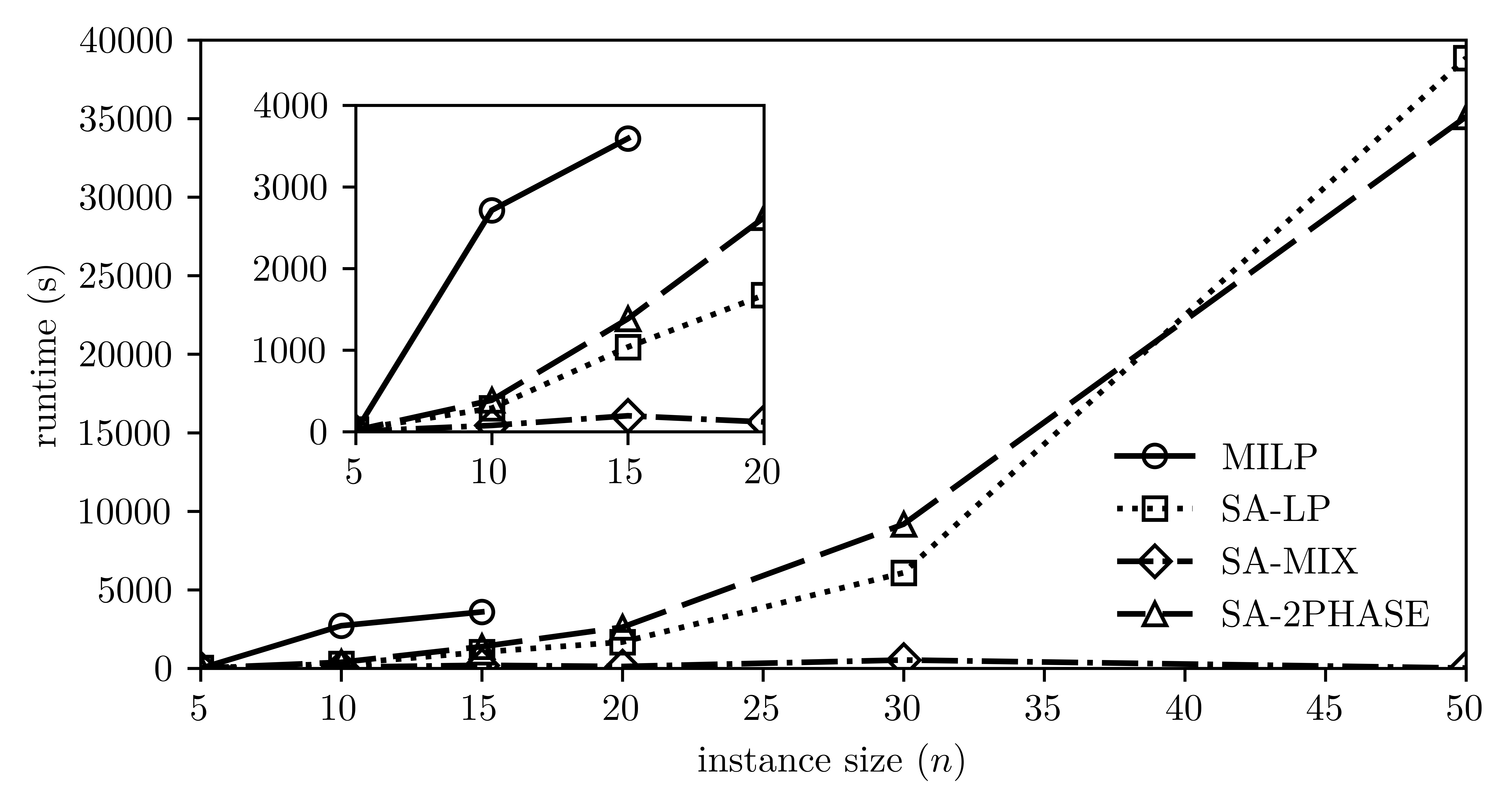}
	\caption{Average runtimes of discussed approaches} \label{fig:gcecsp-runtimes}
\end{figure}

Figure \ref{fig:gcecsp-runtimes} shows the runtime of the approaches. Each point is an average over all instances of a certain size ($n$). For the correct interpretation of the values for the MILP for $n=10,15$, note that the average includes runs that were cut off at 3600 seconds. The MILP was not included for $n \geq 20$ because of very poor performance. Our approach scales better than the MILP, but reports high run times for larger instances as well. However, a reasonable (feasible) solution is often already found relatively quickly, and much time is spent on searching for further improvements. Whenever time is a constraint, cutting of the search after a certain amount of time has elapsed will still yield good results.

\begin{figure}
	\centering
	\includegraphics[width=\textwidth]{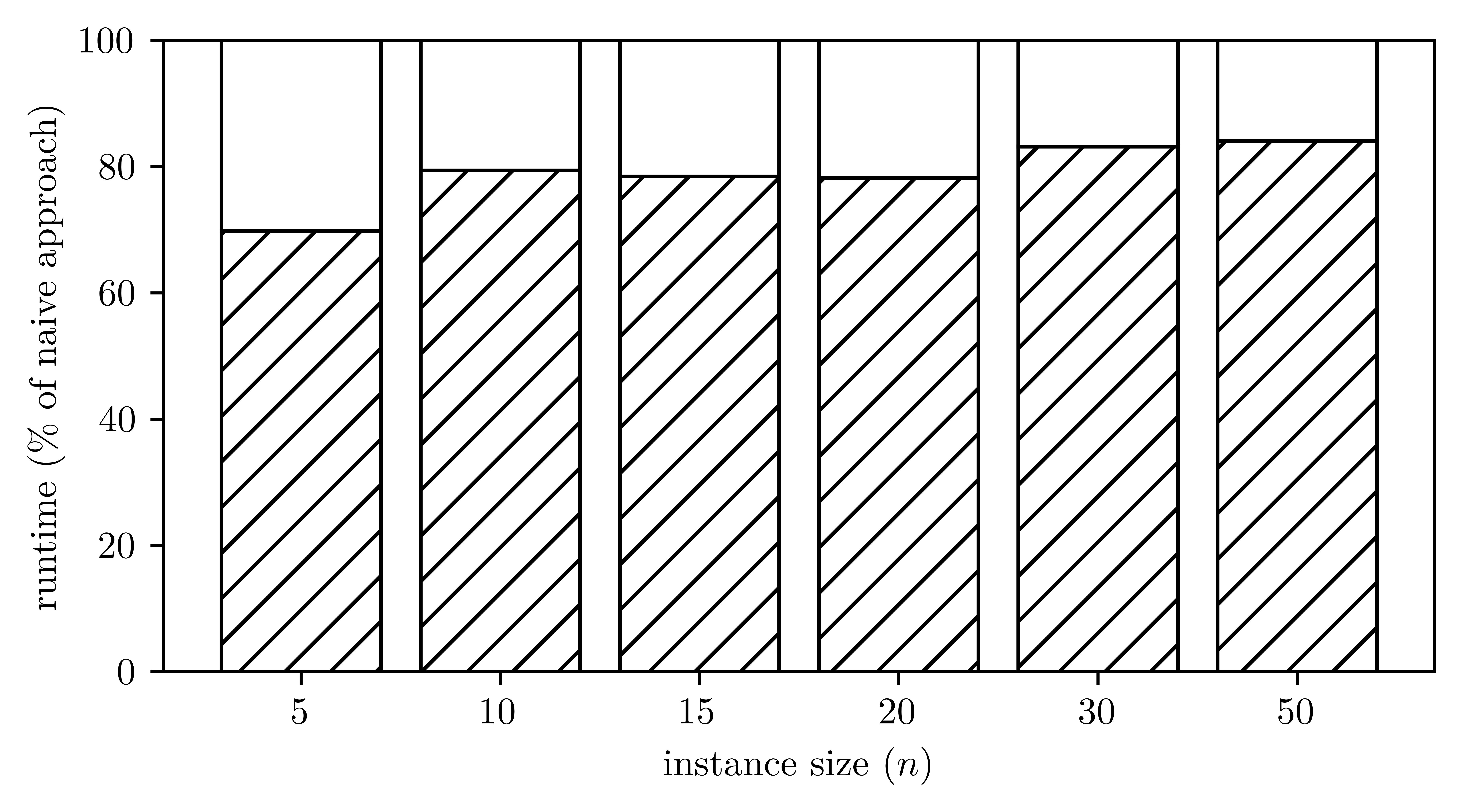}
	\caption{Runtime of the general approach (SA-LP) as a percentage of its naive implementation} \label{fig:gcecsp-naive-comparison}
\end{figure}

The general approach (SA-LP) performs very well on all instance sizes. We are able to reliably find feasible solutions for much larger instances than the MILP. Exchanging the LP in the first part of the search for the `simple' bound estimations (SA-2PHASE) starts to pay off with larger instance sizes. Mixing estimation methods (SA-MIX) does not perform as well as expected. As the iterative update of the LP is much faster than solving the LP from scratch, much of the time gained when using other estimation methods is lost once the LP does have to be solved again. In addition, the quality of the `simple' bounds does seem to be insufficient to guide the search reliably towards good solutions for larger instance sizes. The SA-MIX approach uses only the `simple' lower bound estimation for evaluation whenever these report a positive penalty term. These `simple' bounds do not detect all causes of infeasibility, and are therefore unreliable when the search comes close to a feasible solution. The LP, in contrast, computes the exact penalty term. This may result in moving to candidate solutions with a lower penalty term, but positive `simple' lower bounds, at which point these lower bounds will again be used to evaluate the candidate solutions. We conclude that multiple iterations with the LP are needed to converge on a good solution. It should be noted, that the general approach using only the LP for the evaluation of candidate solutions (SA-LP) also benefits from the fact that the objective value (excluding the penalty term) can be computed before solving the LP. If this is sufficient to reject a candidate solution, the search will move on without solving the LP. To illustrate this point, we ran the general approach (SA-LP) as well as a modified implementation that does not make use of this fact (SA-LP naive) on all instances, seeding the random engine to ensure both approaches follow the exact same path through the search space. Figure~\ref{fig:gcecsp-naive-comparison} shows the average runtime of the general approach (shaded area) as a percentage of the naive approach (full bar). From this, we can see that this simple modification saves around 20\% of computation time.

\begin{figure}
	\centering
	\begin{subfigure}{\textwidth}
		\includegraphics[width=\textwidth]{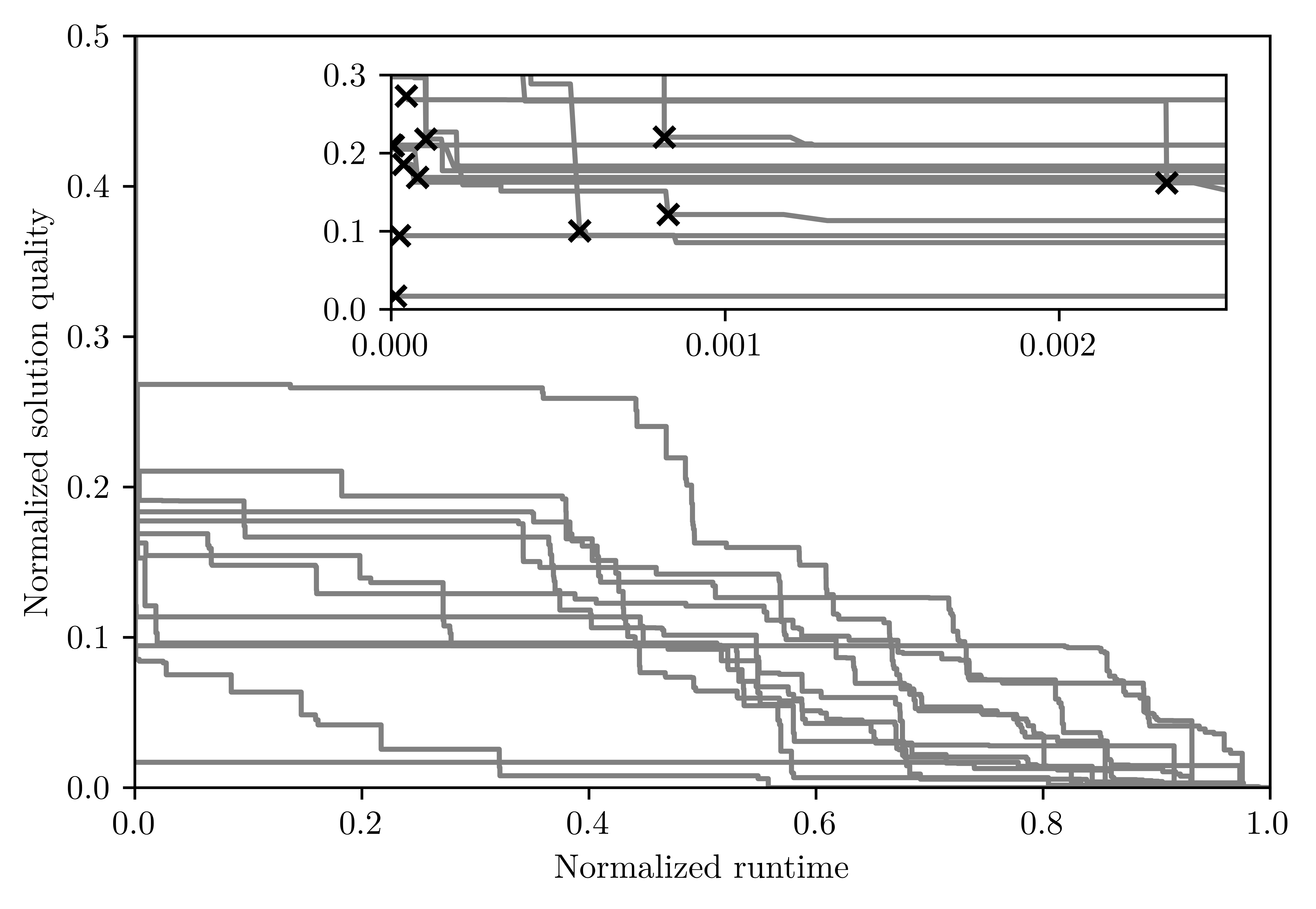}
		\caption{Instances with $n=30$}
		\label{fig:gcecsp-convergence-n30}
	\end{subfigure}
	\begin{subfigure}{\textwidth}
		\includegraphics[width=\textwidth]{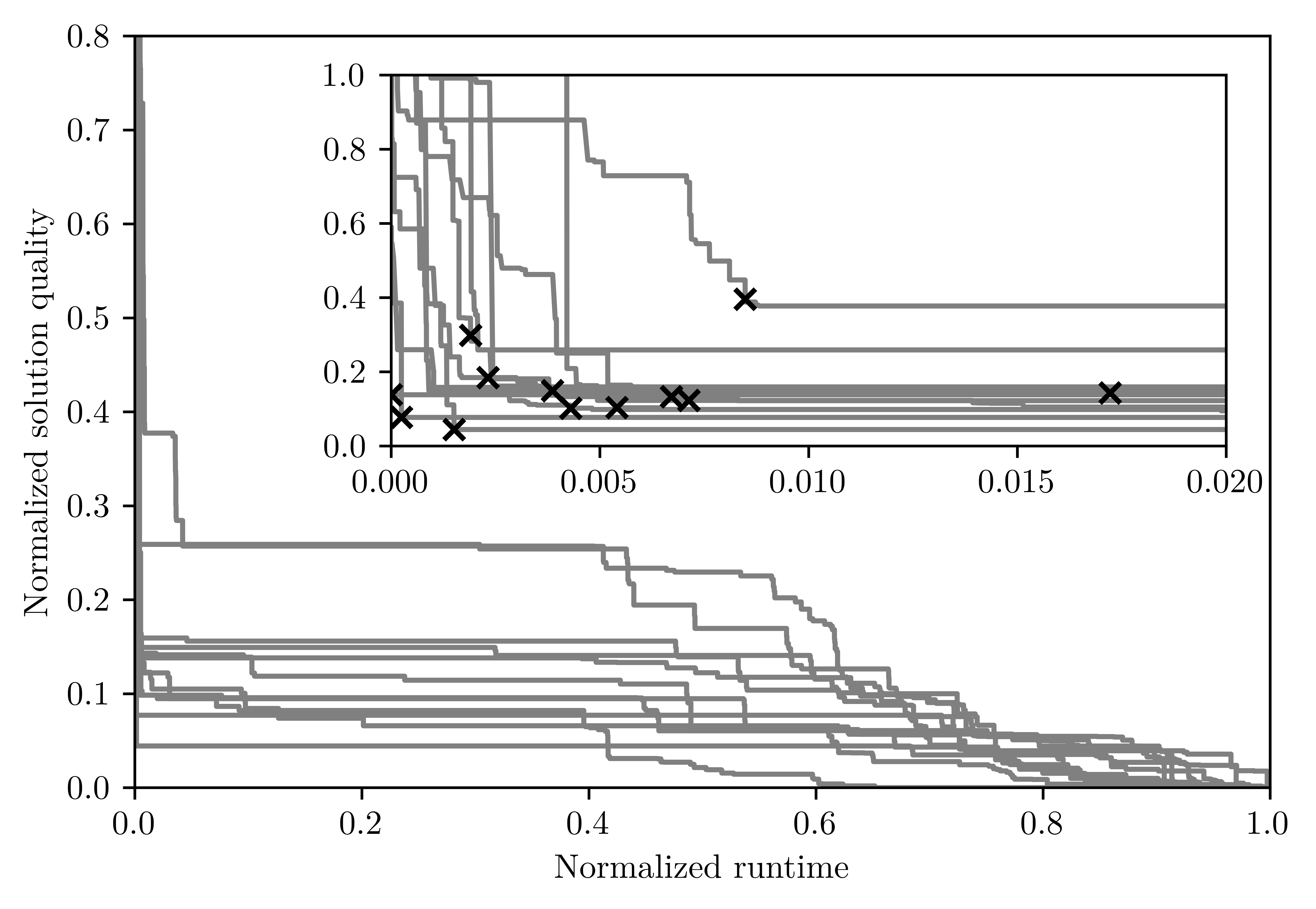}
		\caption{Instances with $n=50$}
		\label{fig:gcecsp-convergence-n50}
	\end{subfigure}
	\caption{Progression of the best found solution over the total runtime}
	\label{fig:gcecsp-convergence}
\end{figure} 

As we mentioned above, cutting of the search early on is a viable approach if time is a constraint. To back up this claim, we reran the SA-LP approach for all instances of size $n=30$ and $n=50$. In Figure \ref{fig:gcecsp-convergence}, we have plotted the progression of the best found solution over the total runtime of the approach. The horizontal axis shows the normalized runtime, and the vertical axis the normalized distance to the best found solution. With the latter, we mean the difference between the current best objective value and the final best objective value expressed as a fraction of the final best objective value. Each run is displayed as a gray line. Note that many lines start far outside the plot, but drop down to within the range displayed after a limited number of successful iterations. The inset zooms in on the first part of the search, and shows the first found feasible solution for each run with a black `x'. From this, we conclude that the first feasible solution is found very quickly, before 0.25\% (for $n=30$) or 2\% (for $n=50$) of the total runtime has passed. We also observe that, after a quick dramatic improvement at the start, the improvement of the best found solution is gradual. Cutting off the search early will provide a good, feasible solution, if time is a constraint.

We already noted in the previous section that finding a good initial solution is important. Even still, the local search will consistently improve the initial solution significantly. Given the sensitivity of the approach to the initial solution, it would be an interesting approach to introduce restarts for instances where the best found solution is not satisfactory. This can be done by disturbing the best found solution using the neighborhood operator a number of times, regardless of the quality of the result, and restarting the search from there.

\section{Conclusions and future work} \label{sec:gcecsp-conclusion}
In this work, we have introduced a hybrid optimization framework for a class of problems that are variants of the GCECSP, which we introduced as a generalization of the CECSP. We have shown the general applicability of our approach, and have studied its performance when applied to the GCECSP with step-wise cost functions. With this, we give an impression of how to tailor the approach to specific variants. Furthermore, we have shown the range of variants that our approach can be applied to, and discussed the required adaptations for a number of cases.

In the application to the variant with step-wise cost functions, we show that we are able to find good solutions for much larger instances, when compared to the MILP. We also observe that it can be used to find reasonable solutions relatively quickly, whereas the MILP starts struggling to find even a feasible solution for instances of size $n=15$ and larger. The fact that we know the part of the objective excluding the penalty term before solving the subproblem, allows us to improve the efficiency of the approach. Exploiting this fact further by replacing the LP for solving the subproblem by alternative approaches has not been as successful.

Furthermore, we open up the possibility of applying the approach to a broader range of problems. Further study of such variants could provide more insight of the quality of the results for specific variants. Further exploration of estimation methods for objectives of Type III (see Section \ref{sec:gcecsp-sub-objectives}) is of particular interest.

\section*{Declarations}

\subsection*{Funding}
No funding was received to assist with the preparation of this manuscript.

\subsection*{Competing interests}
The authors have no competing interests to declare that are relevant to the content of this article.

\subsection*{Data availability}
The instance data used in this paper, as well as an overview of the results presented in this paper are available at: \url{https://doi.org/10.5281/zenodo.10304834}.

\subsection*{Code availability}
The code used to generate the results presented in this paper is available at: \url{https://doi.org/10.5281/zenodo.10304753}.

\begin{appendices}

\section{Full Results}\label{app:gcecsp-results}
This appendix contains the full results of the tests described in Section \ref{sec:gcecsp-results}, in particular those on which the data in Table \ref{tab:gcecsp-results} and Figure \ref{fig:gcecsp-runtimes} is based.

Each table shows the results for all four instances of a particular combination of $n$ and $k$ (displayed in the upper left corner). The first row indicates whether our feasibility test was passed by the instance or not. This is followed by a row indicating the best objective value that we have encountered in any run of the three ($n \geq 20$) or four ($n \leq 15$) listed approaches.

For each of the approaches, we then list the runtime (in seconds), the objective value of the solution found by that particular approach and, in the case of the approaches using simulated annealing, the objective value of the initial solution. In cases where the MILP is cut off because of the time limit set (3600 seconds), this is indicated with the text \texttt{LIMIT}. If no feasible solution is found by the MILP within an hour, there is no objective to report, which is indicated with a `-'. In the row that contains the objective, a value is displayed in \textbf{boldface} if it is the best results among the compared approaches. It is displayed in \textit{italics} if it represents an infeasible solution (i.e., the sum of the violation variables is larger than 0). Values in italics only occur if the approach uses simulated annealing, as the MILP does not include violation variables and therefore does not report infeasible solutions.

Note that the results in these tables are of a single run of each of the approaches, they are not averages or best-of results.
\begin{center}
	\begin{tabular}{ll|rrrr}
		\toprule
		\multicolumn{2}{l|}{$n = 5, k = 2$} & \multicolumn{1}{c}{\# 0} & \multicolumn{1}{c}{\# 1} & \multicolumn{1}{c}{\# 2} & \multicolumn{1}{c}{\# 3} \\
		\midrule
		\midrule
		\multicolumn{2}{l|}{Flow feasible} & \multicolumn{1}{c}{yes} & \multicolumn{1}{c}{yes} & \multicolumn{1}{c}{yes} & \multicolumn{1}{c}{yes} \\
		\midrule
		\multicolumn{2}{l|}{Best known} & 12.42 & 9.74 & 13.08 & 6.37 \\
		\midrule
		\multirow{2}{*}{MILP} & time & 0.43 & 0.44 & 1.03 & 0.67 \\
		& objective & \textbf{12.42} & \textbf{9.74} & \textbf{13.08} & \textbf{6.37} \\
		\midrule
		\multirow{3}{*}{SA-LP} & time & 11.41 & 2.59 & 29.03 & 6.44 \\
		& objective & \textbf{12.42} & 11.66 & \textbf{13.08} & \textbf{6.37} \\
		& init. sol. & 17.24 & 14.58 & 15.50 & 6.37 \\
		\midrule
		\multirow{3}{*}{SA-MIX} & time & 76.10 & 14.43 & 134.66 & 3.25 \\
		& objective & \textbf{12.42} & \textbf{9.74} & \textbf{13.08} & \textbf{6.37} \\
		& init. sol. & 17.24 & 14.58 & 15.50 & 6.37 \\
		\midrule
		\multirow{3}{*}{SA-2PHASE} & time & 2.51 & 3.48 & 13.59 & 1.20 \\
		& objective & \textbf{12.42} & 11.83 & \textbf{13.08} & \textbf{6.37} \\
		& init. sol. & 17.24 & 14.58 & 15.50 & 6.37 \\
		\bottomrule
	\end{tabular}
	
	\begin{tabular}{ll|rrrr}
		\toprule
		\multicolumn{2}{l|}{$n = 5, k = 3$} & \multicolumn{1}{c}{\# 0} & \multicolumn{1}{c}{\# 1} & \multicolumn{1}{c}{\# 2} & \multicolumn{1}{c}{\# 3} \\
		\midrule
		\midrule
		\multicolumn{2}{l|}{Flow feasible} & \multicolumn{1}{c}{yes} & \multicolumn{1}{c}{yes} & \multicolumn{1}{c}{yes} & \multicolumn{1}{c}{yes} \\
		\midrule
		\multicolumn{2}{l|}{Best known} & 15.27 & 12.49 & 17.08 & 19.52 \\
		\midrule
		\multirow{2}{*}{MILP} & time & 0.39 & 0.32 & 0.39 & 0.85 \\
		& objective & \textbf{15.27} & \textbf{12.49} & \textbf{17.08} & \textbf{19.52} \\
		\midrule
		\multirow{3}{*}{SA-LP} & time & 36.37 & 60.17 & 37.33 & 38.98 \\
		& objective & 15.29 & \textbf{12.49} & \textbf{17.08} & 19.53 \\
		& init. sol. & 18.02 & 13.77 & 17.40 & 21.09 \\
		\midrule
		\multirow{3}{*}{SA-MIX} & time & 133.74 & 20.97 & 150.01 & 104.37 \\
		& objective & \textbf{15.27} & \textbf{12.49} & \textbf{17.08} & 19.59 \\
		& init. sol. & 18.02 & 13.77 & 17.40 & 21.09 \\
		\midrule
		\multirow{3}{*}{SA-2PHASE} & time & 52.18 & 23.53 & 28.24 & 32.14 \\
		& objective & \textbf{15.27} & \textbf{12.49} & \textbf{17.08} & 19.53 \\
		& init. sol. & 18.02 & 13.77 & 17.40 & 21.09 \\
		\bottomrule
	\end{tabular}
	
	\begin{tabular}{ll|rrrr}
		\toprule
		\multicolumn{2}{l|}{$n = 5, k = 4$} & \multicolumn{1}{c}{\# 0} & \multicolumn{1}{c}{\# 1} & \multicolumn{1}{c}{\# 2} & \multicolumn{1}{c}{\# 3} \\
		\midrule
		\midrule
		\multicolumn{2}{l|}{Flow feasible} & \multicolumn{1}{c}{yes} & \multicolumn{1}{c}{yes} & \multicolumn{1}{c}{yes} & \multicolumn{1}{c}{yes} \\
		\midrule
		\multicolumn{2}{l|}{Best known} & 14.63 & 16.69 & 9.19 & 10.77 \\
		\midrule
		\multirow{2}{*}{MILP} & time & 0.92 & 12.63 & 0.51 & 0.54 \\
		& objective & \textbf{14.63} & \textbf{16.69} & \textbf{9.19} & \textbf{10.77} \\
		\midrule
		\multirow{3}{*}{SA-LP} & time & 14.20 & 51.16 & 35.30 & 25.74 \\
		& objective & 16.75 & 16.75 & 9.66 & 10.98 \\
		& init. sol. & 167.39 & 19.32 & 10.38 & 17.07 \\
		\midrule
		\multirow{3}{*}{SA-MIX} & time & 0.16 & 106.18 & 208.08 & 28.42 \\
		& objective & \textit{145.89} & 16.75 & \textbf{9.19} & 11.94 \\
		& init. sol. & 167.39 & 19.32 & 10.38 & 17.07 \\
		\midrule
		\multirow{3}{*}{SA-2PHASE} & time & 14.16 & 7.00 & 36.17 & 11.43 \\
		& objective & 15.22 & 18.20 & \textbf{9.19} & 12.80 \\
		& init. sol. & 167.39 & 19.32 & 10.38 & 17.07 \\
		\bottomrule
	\end{tabular}
	
	\begin{tabular}{ll|rrrr}
		\toprule
		\multicolumn{2}{l|}{$n = 10, k = 2$} & \multicolumn{1}{c}{\# 0} & \multicolumn{1}{c}{\# 1} & \multicolumn{1}{c}{\# 2} & \multicolumn{1}{c}{\# 3} \\
		\midrule
		\midrule
		\multicolumn{2}{l|}{Flow feasible} & \multicolumn{1}{c}{yes} & \multicolumn{1}{c}{yes} & \multicolumn{1}{c}{yes} & \multicolumn{1}{c}{yes} \\
		\midrule
		\multicolumn{2}{l|}{Best known} & 19.10 & 17.69 & 15.56 & 33.03 \\
		\midrule
		\multirow{2}{*}{MILP} & time & 3403.32 & \texttt{LIMIT} & 3301.44 & 209.87 \\
		& objective & \textbf{19.10} & \textbf{17.69} & \textbf{15.56} & \textbf{33.03} \\
		\midrule
		\multirow{3}{*}{SA-LP} & time & 177.35 & 124.76 & 245.82 & 184.96 \\
		& objective & 19.43 & 18.20 & \textbf{15.56} & 35.39 \\
		& init. sol. & 243.96 & 308.84 & 18.37 & 38.77 \\
		\midrule
		\multirow{3}{*}{SA-MIX} & time & 171.70 & 0.59 & 199.41 & 123.25 \\
		& objective & 22.15 & \textit{204.66} & \textbf{15.56} & 35.47 \\
		& init. sol. & 243.96 & 308.84 & 18.37 & 38.77 \\
		\midrule
		\multirow{3}{*}{SA-2PHASE} & time & 219.83 & 0.26 & 136.65 & 296.47 \\
		& objective & 21.66 & \textit{85.19} & \textbf{15.56} & 34.45 \\
		& init. sol. & 243.96 & 308.84 & 18.37 & 38.77 \\
		\bottomrule
	\end{tabular}
	
	\begin{tabular}{ll|rrrr}
		\toprule
		\multicolumn{2}{l|}{$n = 10, k = 3$} & \multicolumn{1}{c}{\# 0} & \multicolumn{1}{c}{\# 1} & \multicolumn{1}{c}{\# 2} & \multicolumn{1}{c}{\# 3} \\
		\midrule
		\midrule
		\multicolumn{2}{l|}{Flow feasible} & \multicolumn{1}{c}{yes} & \multicolumn{1}{c}{yes} & \multicolumn{1}{c}{yes} & \multicolumn{1}{c}{yes} \\
		\midrule
		\multicolumn{2}{l|}{Best known} & 21.01 & 28.00 & 26.41 & 22.44 \\
		\midrule
		\multirow{2}{*}{MILP} & time & \texttt{LIMIT} & \texttt{LIMIT} & 402.73 & \texttt{LIMIT} \\
		& objective & \textbf{21.01} & \textbf{28.00} & \textbf{26.41} & \textbf{22.44} \\
		\midrule
		\multirow{3}{*}{SA-LP} & time & 162.92 & 409.86 & 420.81 & 229.25 \\
		& objective & 22.31 & \textbf{28.00} & \textbf{26.41} & 22.80 \\
		& init. sol. & 28.53 & 28.09 & 31.05 & 33.31 \\
		\midrule
		\multirow{3}{*}{SA-MIX} & time & 563.43 & 1270.58 & 1215.83 & 71.37 \\
		& objective & 22.31 & 28.09 & 26.49 & 24.27 \\
		& init. sol. & 28.53 & 28.09 & 31.05 & 33.31 \\
		\midrule
		\multirow{3}{*}{SA-2PHASE} & time & 137.71 & 589.88 & 519.66 & 425.66 \\
		& objective & 22.31 & 28.09 & \textbf{26.41} & 22.80 \\
		& init. sol. & 28.53 & 28.09 & 31.05 & 33.31 \\
		\bottomrule
	\end{tabular}
	
	\begin{tabular}{ll|rrrr}
		\toprule
		\multicolumn{2}{l|}{$n = 10, k = 4$} & \multicolumn{1}{c}{\# 0} & \multicolumn{1}{c}{\# 1} & \multicolumn{1}{c}{\# 2} & \multicolumn{1}{c}{\# 3} \\
		\midrule
		\midrule
		\multicolumn{2}{l|}{Flow feasible} & \multicolumn{1}{c}{yes} & \multicolumn{1}{c}{yes} & \multicolumn{1}{c}{yes} & \multicolumn{1}{c}{yes} \\
		\midrule
		\multicolumn{2}{l|}{Best known} & 27.28 & 25.75 & 29.48 & 24.72 \\
		\midrule
		\multirow{2}{*}{MILP} & time & 41.82 & \texttt{LIMIT} & \texttt{LIMIT} & \texttt{LIMIT} \\
		& objective & \textbf{27.28} & 26.61 & 30.60 & 26.11 \\
		\midrule
		\multirow{3}{*}{SA-LP} & time & 295.05 & 276.10 & 551.23 & 375.96 \\
		& objective & 27.32 & 26.49 & \textbf{29.48} & \textbf{24.72} \\
		& init. sol. & 30.70 & 31.10 & 114.09 & 203.20 \\
		\midrule
		\multirow{3}{*}{SA-MIX} & time & 13.84 & 29.44 & 209.57 & 15.33 \\
		& objective & 30.70 & 28.47 & 32.70 & \textit{31.82} \\
		& init. sol. & 30.70 & 31.10 & 114.09 & 203.20 \\
		\midrule
		\multirow{3}{*}{SA-2PHASE} & time & 545.23 & 549.70 & 597.15 & 619.47 \\
		& objective & 27.32 & \textbf{25.75} & 29.94 & 25.41 \\
		& init. sol. & 30.70 & 31.10 & 114.09 & 203.20 \\
		\bottomrule
	\end{tabular}
	
	\begin{tabular}{ll|rrrr}
		\toprule
		\multicolumn{2}{l|}{$n = 15, k = 2$} & \multicolumn{1}{c}{\# 0} & \multicolumn{1}{c}{\# 1} & \multicolumn{1}{c}{\# 2} & \multicolumn{1}{c}{\# 3} \\
		\midrule
		\midrule
		\multicolumn{2}{l|}{Flow feasible} & \multicolumn{1}{c}{yes} & \multicolumn{1}{c}{yes} & \multicolumn{1}{c}{yes} & \multicolumn{1}{c}{yes} \\
		\midrule
		\multicolumn{2}{l|}{Best known} & 45.32 & 29.64 & 32.52 & 42.50 \\
		\midrule
		\multirow{2}{*}{MILP} & time & 3472.54 & \texttt{LIMIT} & \texttt{LIMIT} & \texttt{LIMIT} \\
		& objective & \textbf{45.32} & 29.68 & \textit{-} & 43.18 \\
		\midrule
		\multirow{3}{*}{SA-LP} & time & 1166.54 & 826.58 & 694.15 & 756.64 \\
		& objective & \textbf{45.32} & \textbf{29.64} & \textbf{32.52} & \textbf{42.50} \\
		& init. sol. & 89.63 & 37.38 & 147.25 & 51.33 \\
		\midrule
		\multirow{3}{*}{SA-MIX} & time & 141.57 & 94.75 & 0.29 & 200.75 \\
		& objective & 45.43 & 36.14 & 48.77 & 50.56 \\
		& init. sol. & 89.63 & 37.38 & 147.25 & 51.33 \\
		\midrule
		\multirow{3}{*}{SA-2PHASE} & time & 1954.29 & 1333.77 & 432.17 & 486.36 \\
		& objective & \textbf{45.32} & \textbf{29.64} & 33.02 & 45.61 \\
		& init. sol. & 89.63 & 37.38 & 147.25 & 51.33 \\
		\bottomrule
	\end{tabular}
	
	\begin{tabular}{ll|rrrr}
		\toprule
		\multicolumn{2}{l|}{$n = 15, k = 3$} & \multicolumn{1}{c}{\# 0} & \multicolumn{1}{c}{\# 1} & \multicolumn{1}{c}{\# 2} & \multicolumn{1}{c}{\# 3} \\
		\midrule
		\midrule
		\multicolumn{2}{l|}{Flow feasible} & \multicolumn{1}{c}{yes} & \multicolumn{1}{c}{yes} & \multicolumn{1}{c}{yes} & \multicolumn{1}{c}{yes} \\
		\midrule
		\multicolumn{2}{l|}{Best known} & 37.18 & 40.41 & 45.08 & 32.17 \\
		\midrule
		\multirow{2}{*}{MILP} & time & \texttt{LIMIT} & \texttt{LIMIT} & \texttt{LIMIT} & \texttt{LIMIT} \\
		& objective & 39.03 & \textit{-} & \textit{-} & 38.56 \\
		\midrule
		\multirow{3}{*}{SA-LP} & time & 837.38 & 994.47 & 1295.04 & 1293.10 \\
		& objective & 37.56 & \textbf{40.41} & 45.31 & 32.83 \\
		& init. sol. & 50.79 & 50.00 & 105.65 & 40.50 \\
		\midrule
		\multirow{3}{*}{SA-MIX} & time & 942.28 & 63.18 & 23.69 & 41.43 \\
		& objective & 40.95 & 45.98 & 49.04 & 40.50 \\
		& init. sol. & 50.79 & 50.00 & 105.65 & 40.50 \\
		\midrule
		\multirow{3}{*}{SA-2PHASE} & time & 817.44 & 969.81 & 1671.45 & 2317.23 \\
		& objective & \textbf{37.18} & 42.93 & \textbf{45.08} & \textbf{32.17} \\
		& init. sol. & 50.79 & 50.00 & 105.65 & 40.50 \\
		\bottomrule
	\end{tabular}
	
	\begin{tabular}{ll|rrrr}
		\toprule
		\multicolumn{2}{l|}{$n = 15, k = 4$} & \multicolumn{1}{c}{\# 0} & \multicolumn{1}{c}{\# 1} & \multicolumn{1}{c}{\# 2} & \multicolumn{1}{c}{\# 3} \\
		\midrule
		\midrule
		\multicolumn{2}{l|}{Flow feasible} & \multicolumn{1}{c}{yes} & \multicolumn{1}{c}{yes} & \multicolumn{1}{c}{yes} & \multicolumn{1}{c}{yes} \\
		\midrule
		\multicolumn{2}{l|}{Best known} & 42.38 & 41.76 & 29.34 & 33.96 \\
		\midrule
		\multirow{2}{*}{MILP} & time & \texttt{LIMIT} & \texttt{LIMIT} & \texttt{LIMIT} & \texttt{LIMIT} \\
		& objective & 42.68 & \textit{-} & 31.57 & \textit{-} \\
		\midrule
		\multirow{3}{*}{SA-LP} & time & 1280.55 & 1504.55 & 1017.33 & 787.37 \\
		& objective & \textbf{42.38} & \textbf{41.76} & 29.61 & \textbf{33.96} \\
		& init. sol. & 51.96 & 54.36 & 82.57 & 42.93 \\
		\midrule
		\multirow{3}{*}{SA-MIX} & time & 17.26 & 129.05 & 0.99 & 46.33 \\
		& objective & 51.96 & 50.36 & \textit{68.28} & 42.93 \\
		& init. sol. & 51.96 & 54.36 & 82.57 & 42.93 \\
		\midrule
		\multirow{3}{*}{SA-2PHASE} & time & 1419.34 & 1909.92 & 2076.76 & 1216.30 \\
		& objective & 43.84 & 44.17 & \textbf{29.34} & 34.58 \\
		& init. sol. & 51.96 & 54.36 & 82.57 & 42.93 \\
		\bottomrule
	\end{tabular}
	
	\begin{tabular}{ll|rrrr}
		\toprule
		\multicolumn{2}{l|}{$n = 20, k = 2$} & \multicolumn{1}{c}{\# 0} & \multicolumn{1}{c}{\# 1} & \multicolumn{1}{c}{\# 2} & \multicolumn{1}{c}{\# 3} \\
		\midrule
		\midrule
		\multicolumn{2}{l|}{Flow feasible} & \multicolumn{1}{c}{yes} & \multicolumn{1}{c}{yes} & \multicolumn{1}{c}{yes} & \multicolumn{1}{c}{yes} \\
		\midrule
		\multicolumn{2}{l|}{Best known} & 52.47 & 49.03 & 49.14 & 44.29 \\
		\midrule
		\multirow{3}{*}{SA-LP} & time & 1502.61 & 1708.22 & 1331.91 & 1156.45 \\
		& objective & 54.18 & 49.11 & 50.00 & \textbf{44.29} \\
		& init. sol. & 195.25 & 319.42 & 59.13 & 278.03 \\
		\midrule
		\multirow{3}{*}{SA-MIX} & time & 13.42 & 2.71 & 601.63 & 155.35 \\
		& objective & \textit{62.79} & \textit{318.19} & \textit{55.74} & 46.54 \\
		& init. sol. & 195.25 & 319.42 & 59.13 & 278.03 \\
		\midrule
		\multirow{3}{*}{SA-2PHASE} & time & 3198.27 & 3191.63 & 1373.58 & 1332.58 \\
		& objective & \textbf{52.47} & \textbf{49.03} & \textbf{49.14} & 44.70 \\
		& init. sol. & 195.25 & 319.42 & 59.13 & 278.03 \\
		\bottomrule
	\end{tabular}
	
	\begin{tabular}{ll|rrrr}
		\toprule
		\multicolumn{2}{l|}{$n = 20, k = 3$} & \multicolumn{1}{c}{\# 0} & \multicolumn{1}{c}{\# 1} & \multicolumn{1}{c}{\# 2} & \multicolumn{1}{c}{\# 3} \\
		\midrule
		\midrule
		\multicolumn{2}{l|}{Flow feasible} & \multicolumn{1}{c}{yes} & \multicolumn{1}{c}{yes} & \multicolumn{1}{c}{yes} & \multicolumn{1}{c}{yes} \\
		\midrule
		\multicolumn{2}{l|}{Best known} & 65.33 & 57.83 & 52.91 & 62.42 \\
		\midrule
		\multirow{3}{*}{SA-LP} & time & 1542.54 & 1463.81 & 1797.87 & 1969.42 \\
		& objective & 65.92 & 58.01 & \textbf{52.91} & \textbf{62.42} \\
		& init. sol. & 192.59 & 170.00 & 161.81 & 178.35 \\
		\midrule
		\multirow{3}{*}{SA-MIX} & time & 80.95 & 4.26 & 131.92 & 7.42 \\
		& objective & 68.99 & \textit{149.90} & 63.37 & \textit{113.21} \\
		& init. sol. & 192.59 & 170.00 & 161.81 & 178.35 \\
		\midrule
		\multirow{3}{*}{SA-2PHASE} & time & 1726.78 & 2661.88 & 2985.21 & 3246.13 \\
		& objective & \textbf{65.33} & \textbf{57.83} & 53.54 & 62.51 \\
		& init. sol. & 192.59 & 170.00 & 161.81 & 178.35 \\
		\bottomrule
	\end{tabular}
	
	\begin{tabular}{ll|rrrr}
		\toprule
		\multicolumn{2}{l|}{$n = 20, k = 4$} & \multicolumn{1}{c}{\# 0} & \multicolumn{1}{c}{\# 1} & \multicolumn{1}{c}{\# 2} & \multicolumn{1}{c}{\# 3} \\
		\midrule
		\midrule
		\multicolumn{2}{l|}{Flow feasible} & \multicolumn{1}{c}{yes} & \multicolumn{1}{c}{yes} & \multicolumn{1}{c}{yes} & \multicolumn{1}{c}{yes} \\
		\midrule
		\multicolumn{2}{l|}{Best known} & 62.09 & 45.71 & 53.69 & 58.47 \\
		\midrule
		\multirow{3}{*}{SA-LP} & time & 2162.59 & 1909.21 & 1930.20 & 1621.52 \\
		& objective & 63.22 & \textbf{45.71} & 55.54 & 58.54 \\
		& init. sol. & 165.04 & 105.35 & 314.93 & 70.27 \\
		\midrule
		\multirow{3}{*}{SA-MIX} & time & 2.84 & 1.34 & 1.89 & 36.27 \\
		& objective & \textit{146.67} & \textit{77.84} & \textit{250.40} & 66.43 \\
		& init. sol. & 165.04 & 105.35 & 314.93 & 70.27 \\
		\midrule
		\multirow{3}{*}{SA-2PHASE} & time & 2542.63 & 2794.88 & 3813.35 & 2615.45 \\
		& objective & \textbf{62.09} & 45.85 & \textbf{53.69} & \textbf{58.47} \\
		& init. sol. & 165.04 & 105.35 & 314.93 & 70.27 \\
		\bottomrule
	\end{tabular}
	
	\begin{tabular}{ll|rrrr}
		\toprule
		\multicolumn{2}{l|}{$n = 30, k = 2$} & \multicolumn{1}{c}{\# 0} & \multicolumn{1}{c}{\# 1} & \multicolumn{1}{c}{\# 2} & \multicolumn{1}{c}{\# 3} \\
		\midrule
		\midrule
		\multicolumn{2}{l|}{Flow feasible} & \multicolumn{1}{c}{yes} & \multicolumn{1}{c}{yes} & \multicolumn{1}{c}{yes} & \multicolumn{1}{c}{yes} \\
		\midrule
		\multicolumn{2}{l|}{Best known} & 79.77 & 78.99 & 75.59 & 69.29 \\
		\midrule
		\multirow{3}{*}{SA-LP} & time & 5987.90 & 5293.30 & 7501.06 & 4438.27 \\
		& objective & \textbf{79.77} & 81.79 & \textbf{75.59} & 70.03 \\
		& init. sol. & 201.16 & 101.66 & 79.57 & 106.25 \\
		\midrule
		\multirow{3}{*}{SA-MIX} & time & 79.23 & 492.44 & 282.40 & 81.74 \\
		& objective & 89.98 & 82.82 & 79.57 & 77.21 \\
		& init. sol. & 201.16 & 101.66 & 79.57 & 106.25 \\
		\midrule
		\multirow{3}{*}{SA-2PHASE} & time & 9443.22 & 8515.09 & 7840.55 & 8913.61 \\
		& objective & 83.26 & \textbf{78.99} & 78.20 & \textbf{69.29} \\
		& init. sol. & 201.16 & 101.66 & 79.57 & 106.25 \\
		\bottomrule
	\end{tabular}
	
	\begin{tabular}{ll|rrrr}
		\toprule
		\multicolumn{2}{l|}{$n = 30, k = 3$} & \multicolumn{1}{c}{\# 0} & \multicolumn{1}{c}{\# 1} & \multicolumn{1}{c}{\# 2} & \multicolumn{1}{c}{\# 3} \\
		\midrule
		\midrule
		\multicolumn{2}{l|}{Flow feasible} & \multicolumn{1}{c}{yes} & \multicolumn{1}{c}{yes} & \multicolumn{1}{c}{yes} & \multicolumn{1}{c}{yes} \\
		\midrule
		\multicolumn{2}{l|}{Best known} & 81.40 & 85.20 & 74.27 & 85.25 \\
		\midrule
		\multirow{3}{*}{SA-LP} & time & 6495.80 & 5113.45 & 6655.20 & 6148.46 \\
		& objective & \textbf{81.40} & 86.49 & \textbf{74.27} & 85.64 \\
		& init. sol. & 360.08 & 100.75 & 91.16 & 96.37 \\
		\midrule
		\multirow{3}{*}{SA-MIX} & time & 4.52 & 155.68 & 103.69 & 26.29 \\
		& objective & \textit{130.42} & 94.67 & 91.16 & 92.82 \\
		& init. sol. & 360.08 & 100.75 & 91.16 & 96.37 \\
		\midrule
		\multirow{3}{*}{SA-2PHASE} & time & 6833.45 & 9203.74 & 9654.89 & 8082.87 \\
		& objective & 83.55 & \textbf{85.20} & 77.64 & \textbf{85.25} \\
		& init. sol. & 360.08 & 100.75 & 91.16 & 96.37 \\
		\bottomrule
	\end{tabular}
	
	\begin{tabular}{ll|rrrr}
		\toprule
		\multicolumn{2}{l|}{$n = 30, k = 4$} & \multicolumn{1}{c}{\# 0} & \multicolumn{1}{c}{\# 1} & \multicolumn{1}{c}{\# 2} & \multicolumn{1}{c}{\# 3} \\
		\midrule
		\midrule
		\multicolumn{2}{l|}{Flow feasible} & \multicolumn{1}{c}{yes} & \multicolumn{1}{c}{yes} & \multicolumn{1}{c}{yes} & \multicolumn{1}{c}{yes} \\
		\midrule
		\multicolumn{2}{l|}{Best known} & 75.34 & 75.56 & 83.54 & 73.18 \\
		\midrule
		\multirow{3}{*}{SA-LP} & time & 5793.06 & 8548.83 & 4606.60 & 6463.71 \\
		& objective & 75.64 & \textbf{75.56} & 86.28 & 73.46 \\
		& init. sol. & 94.59 & 92.29 & 143.81 & 83.83 \\
		\midrule
		\multirow{3}{*}{SA-MIX} & time & 123.32 & 264.12 & 114.96 & 104.98 \\
		& objective & 92.28 & \textit{89.89} & 105.69 & 83.83 \\
		& init. sol. & 94.59 & 92.29 & 143.81 & 83.83 \\
		\midrule
		\multirow{3}{*}{SA-2PHASE} & time & 8152.17 & 14057.77 & 7972.45 & 11481.82 \\
		& objective & \textbf{75.34} & 76.30 & \textbf{83.54} & \textbf{73.18} \\
		& init. sol. & 94.59 & 92.29 & 143.81 & 83.83 \\
		\bottomrule
	\end{tabular}
	
	\begin{tabular}{ll|rrrr}
		\toprule
		\multicolumn{2}{l|}{$n = 50, k = 2$} & \multicolumn{1}{c}{\# 0} & \multicolumn{1}{c}{\# 1} & \multicolumn{1}{c}{\# 2} & \multicolumn{1}{c}{\# 3} \\
		\midrule
		\midrule
		\multicolumn{2}{l|}{Flow feasible} & \multicolumn{1}{c}{yes} & \multicolumn{1}{c}{yes} & \multicolumn{1}{c}{yes} & \multicolumn{1}{c}{yes} \\
		\midrule
		\multicolumn{2}{l|}{Best known} & 128.78 & 118.32 & 121.08 & 127.03 \\
		\midrule
		\multirow{3}{*}{SA-LP} & time & 41080.91 & 38014.12 & 47866.84 & 30884.50 \\
		& objective & \textbf{128.78} & \textbf{118.32} & \textbf{121.08} & \textbf{127.03} \\
		& init. sol. & 383.00 & 255.24 & 2634.19 & 180.49 \\
		\midrule
		\multirow{3}{*}{SA-MIX} & time & 38.09 & 668.13 & 6.05 & 288.01 \\
		& objective & \textit{362.09} & \textit{130.84} & \textit{2307.93} & 135.33 \\
		& init. sol. & 383.00 & 255.24 & 2634.19 & 180.49 \\
		\midrule
		\multirow{3}{*}{SA-2PHASE} & time & 39614.37 & 35312.87 & 34989.71 & 34824.88 \\
		& objective & 131.76 & 120.82 & 122.73 & 128.02 \\
		& init. sol. & 383.00 & 255.24 & 2634.19 & 180.49 \\
		\bottomrule
	\end{tabular}
	
	\begin{tabular}{ll|rrrr}
		\toprule
		\multicolumn{2}{l|}{$n = 50, k = 3$} & \multicolumn{1}{c}{\# 0} & \multicolumn{1}{c}{\# 1} & \multicolumn{1}{c}{\# 2} & \multicolumn{1}{c}{\# 3} \\
		\midrule
		\midrule
		\multicolumn{2}{l|}{Flow feasible} & \multicolumn{1}{c}{yes} & \multicolumn{1}{c}{yes} & \multicolumn{1}{c}{yes} & \multicolumn{1}{c}{yes} \\
		\midrule
		\multicolumn{2}{l|}{Best known} & 134.60 & 141.69 & 139.51 & 152.90 \\
		\midrule
		\multirow{3}{*}{SA-LP} & time & 44735.46 & 42392.18 & 34161.03 & 25714.16 \\
		& objective & 135.03 & \textbf{141.69} & 140.43 & 154.21 \\
		& init. sol. & 677.57 & 161.94 & 661.56 & 247.66 \\
		\midrule
		\multirow{3}{*}{SA-MIX} & time & 2.87 & 524.66 & 9.37 & 174.86 \\
		& objective & \textit{475.18} & 161.94 & \textit{278.72} & \textit{205.48} \\
		& init. sol. & 677.57 & 161.94 & 661.56 & 247.66 \\
		\midrule
		\multirow{3}{*}{SA-2PHASE} & time & 31956.34 & 37477.15 & 34576.79 & 26084.72 \\
		& objective & \textbf{134.60} & 145.87 & \textbf{139.51} & \textbf{152.90} \\
		& init. sol. & 677.57 & 161.94 & 661.56 & 247.66 \\
		\bottomrule
	\end{tabular}
	
	\begin{tabular}{ll|rrrr}
		\toprule
		\multicolumn{2}{l|}{$n = 50, k = 4$} & \multicolumn{1}{c}{\# 0} & \multicolumn{1}{c}{\# 1} & \multicolumn{1}{c}{\# 2} & \multicolumn{1}{c}{\# 3} \\
		\midrule
		\midrule
		\multicolumn{2}{l|}{Flow feasible} & \multicolumn{1}{c}{yes} & \multicolumn{1}{c}{yes} & \multicolumn{1}{c}{yes} & \multicolumn{1}{c}{yes} \\
		\midrule
		\multicolumn{2}{l|}{Best known} & 133.09 & 121.30 & 134.82 & 125.22 \\
		\midrule
		\multirow{3}{*}{SA-LP} & time & 42071.77 & 33440.85 & 38900.98 & 46989.68 \\
		& objective & 135.06 & \textbf{121.30} & \textbf{134.82} & 131.50 \\
		& init. sol. & 501.21 & 292.40 & 236.23 & 362.49 \\
		\midrule
		\multirow{3}{*}{SA-MIX} & time & 7.28 & 201.64 & 111.17 & 62.61 \\
		& objective & \textit{480.72} & \textit{161.28} & \textit{211.58} & \textit{362.49} \\
		& init. sol. & 501.21 & 292.40 & 236.23 & 362.49 \\
		\midrule
		\multirow{3}{*}{SA-2PHASE} & time & 41332.60 & 27681.98 & 39418.33 & 38166.36 \\
		& objective & \textbf{133.09} & 121.89 & 138.18 & \textbf{125.22} \\
		& init. sol. & 501.21 & 292.40 & 236.23 & 362.49 \\
		\bottomrule
	\end{tabular}
\end{center}

\end{appendices}

\printbibliography

\end{document}